\documentclass[12pt,a4paper,reqno]{amsart}
\usepackage[T1]{fontenc}

\usepackage{lmodern}
\usepackage[a4paper,textwidth=16.5cm,textheight=23cm,centering]{geometry}
\usepackage{amsmath,amsfonts,amssymb,amsthm,mathtools,amsrefs}
\usepackage[shortlabels]{enumitem}
\usepackage{graphicx,booktabs,array,needspace}
\usepackage[expansion=false]{microtype}
\usepackage{xcolor}
\definecolor{CitationBlue}{RGB}{0,0,192}
\definecolor{CrossRefWine}{RGB}{128,0,32}
\usepackage[colorlinks=true,citecolor=CitationBlue,linkcolor=CrossRefWine,urlcolor=black]{hyperref}
\hypersetup{pdftitle={Vanishing transport noise enhances stability of stratified Couette flow in the 3D Boussinesq system},pdfauthor={Lei Zhang},pdfsubject={Projected transport noise, enhanced dissipation, and nonlinear buoyancy feedback}}
\setlist{itemsep=5pt,topsep=6pt,parsep=2pt}
\allowdisplaybreaks[2]
\numberwithin{equation}{section}
\newtheorem{theorem}{Theorem}[section]
\newtheorem{proposition}[theorem]{Proposition}
\newtheorem{lemma}[theorem]{Lemma}
\newtheorem{corollary}[theorem]{Corollary}
\theoremstyle{definition}

\newtheorem{remark}[theorem]{Remark}
\newtheorem*{remark*}{Remark}

\newcommand{\R}{\mathbb R}
\newcommand{\T}{\mathbb T}
\newcommand{\E}{\mathbb E}
\newcommand{\Prob}{\mathbb P}
\newcommand{\cF}{\mathcal F}
\newcommand{\cA}{\mathcal A}
\newcommand{\ad}{\mathrm{ad}}
\newcommand{\dd}{\mathrm d}
\newcommand{\eps}{\varepsilon}
\newcommand{\av}[1]{\langle #1\rangle}
\newcommand{\norm}[1]{\left\|#1\right\|}
\newcommand{\sh}{\sharp}
\DeclareMathOperator{\diver}{div}
\let\cal\mathcal

\title[Stratified Couette flow with transport noise]{Vanishing transport noise enhances stability of stratified Couette flow in the 3D Boussinesq system}
\author{Lei Zhang}
\address{School of Mathematics and Statistics, Hubei Key Laboratory of Engineering Modeling and Scientific Computing, Huazhong University of Science and Technology}
\email{lei\_zhang@hust.edu.cn}
\date{}

\begin{document}
\begin{abstract}
In this paper, we study the enhanced nonlinear stability of stably
stratified Couette flow in the three-dimensional Boussinesq system driven
by a common affine Stratonovich transport noise. Our main result shows
that, for fixed $\beta>1/2$, suitably chosen noise amplitudes enlarge the
range of perturbations for which we prove stability, even as they vanish
with viscosity. For each fixed deterministic Sobolev perturbation
satisfying the sufficient smallness conditions, the strong solution is
global and obeys stability estimates on an event of prescribed high
probability. No spatial localization or vanishing initial double means
are required. The proof combines noise-induced mixing with coupled energy
estimates to control nonlinear feedback and preserve the physical damping
of the fast vertical velocity and temperature.
\end{abstract}

\maketitle

\section{Introduction}

\subsection{Formulation of the problem}
In this paper, we study how vanishing transport noise improves sufficient
conditions for nonlinear stability of stratified Couette flow in the
three-dimensional Boussinesq system on
$\Omega=\mathbb T_x\times\mathbb R_y\times\mathbb T_z$, with both periodic
directions of period $2\pi$. To formulate the problem, we start from the
classical equations and model unresolved fluid motion through a common
random advection.

Let $y$ denote the upward vertical coordinate, and let $U$, $b$ and $\Pi$ denote the total velocity, upward buoyancy acceleration and kinematic pressure, respectively. We define $e_x=(1,0,0)$, $e_y=(0,1,0)$, and $e_z=(0,0,1)$. For thermal stratification, $b$ is proportional to the temperature deviation from a reference state. The classical Boussinesq
equations are given by
\[
 \left\{\begin{aligned}
 \partial_tU+(U\cdot\nabla)U+\nabla\Pi
       &=\nu_{\rm phys}\Delta U+b e_y,\\
 \partial_tb+U\cdot\nabla b&=\kappa_{\rm phys}\Delta b,\\
 \nabla\cdot U&=0.
 \end{aligned}\right.
\]
Let us take a background shear rate $S>0$ and a buoyancy frequency $N>0$.
The constant background buoyancy gradient $N^2$ is positive, so that
the stratification is stable. Using length scale $L$
and time scale $S^{-1}$, we scale velocity, buoyancy and pressure by
$SL$, $S^2L$ and $S^2L^2$, and retain $U,b,\Pi$ for the dimensionless
fields. We assume $\kappa_{\rm phys}=\nu_{\rm phys}$ (Prandtl number one)
and set \mbox{$\nu=\nu_{\rm phys}/(SL^2)>0$} and $\beta=N/S>0$.
On $\Omega$,  the background profiles are
$\overline U=y e_x$, $\overline b=\beta^2y$ and
$\overline\Pi=\frac12\beta^2y^2$, which satisfy hydrostatic balance and
have gradient Richardson number $Ri=\beta^2$.

Unresolved motions also advect the resolved velocity and buoyancy fields.
To represent this effect, we introduce random particle displacements
as in the location-uncertainty formulations of M\'emin~\cite{Memin2014}
and Resseguier et al.~\cite{ResseguierMeminChapron2017}.
Let $W$ be a real Brownian motion on a filtered probability space
satisfying the usual conditions, and let $\circ$ denote Stratonovich
integration. We consider the particle motion
\[
 \dd\mathbf X_t=U(\mathbf X_t,t)\dd t
       +\sigma(\mathbf X_t)\circ\dd W_t,
 \quad \sigma(x,y,z)=\rho y e_z,\quad \rho\ge0.
\]
Here the divergence-free field $\sigma$ describes a random spanwise shear,  and the
Stratonovich integration preserves the chain rule along these
trajectories.  We refer to the work by Flandoli and
Pappalettera~\cite{FlandoliPappalettera2021} for a scale-separation justification of transport noise in a
coupled two-dimensional hydrodynamic system.

\Needspace{7\baselineskip}
Using the same random displacement to transport velocity components
and buoyancy, we adopt the componentwise transport closure
\[
 \left\{\begin{aligned}
 \dd U+[(U\cdot\nabla)U+\nabla\Pi]\dd t
       +[(\sigma\cdot\nabla)U+\nabla\chi]\circ\dd W_t
       &=(\nu\Delta U+b e_y)\dd t,\\
 \dd b+(U\cdot\nabla b)\dd t
       +(\sigma\cdot\nabla b)\circ\dd W_t
       &=\nu\Delta b\dd t,\\
 \nabla\cdot U&=0,
 \end{aligned}\right.
\]
where the stochastic pressure coefficient $\chi$ enforces
incompressibility. Since
$\sigma\cdot\nabla\overline U=\sigma\cdot\nabla\overline b=0$,
the Couette flow and linear stratification remain exact background
solutions under this common random transport.
Set $U=\overline U+u$, $b=\overline b-\beta\theta$ and
$p=\Pi-\overline\Pi$, where \mbox{$u=(u^1,u^2,u^3)$}.
Expanding the velocity advection gives
$$(U\cdot\nabla)U=y\partial_xu+u^2e_x+(u\cdot\nabla)u.$$
The terms $y\partial_xu$ and $u^2e_x$ describe background advection
and the effect of vertical motion across the shear, respectively.
The buoyancy equation produces the scalar source $+\beta u^2$,
while the restoring force in the velocity equation is
$-\beta\theta e_y$. These terms cancel in the combined
velocity--temperature $L^2$ energy. Incompressibility also gives
$(u\cdot\nabla)u=\diver(u\otimes u)$ and
$(u\cdot\nabla)\theta=\diver(u\theta)$.

After subtracting the background equations, we apply the Leray
projection $\mathcal P$ on $\Omega$ to remove the perturbation pressures.
Note that although $\sigma$ is divergence-free, the term
\mbox{$\diver(y\partial_zu)=\partial_zu^2$} need not vanish.
The stochastic pressure therefore satisfies
$\Delta\chi=-\rho\partial_zu^2$, which yields
$\rho y\partial_zu+\nabla\chi=\rho\mathcal P(y\partial_zu)$.
For pure Leray-projected transport noise, see  for instance Goodair~\cite{Goodair2025}.

We thus obtain the stochastic perturbation system
\begin{equation}\label{eq:model}
\left\{\begin{aligned}
 \dd u+\mathcal P\big[y\partial_xu+u^2e_x
             +(u\cdot\nabla)u+\beta\theta e_y\big]\dd t
       +\rho\mathcal P(y\partial_zu)\circ\dd W_t
       &=\nu\Delta u\dd t,\\
 \dd\theta+\big[y\partial_x\theta+(u\cdot\nabla)\theta
                    -\beta u^2\big]\dd t
       +\rho y\partial_z\theta\circ\dd W_t
       &=\nu\Delta\theta\dd t,\\
 \nabla\cdot u&=0,\\
 (u,\theta)|_{t=0}&=(u_{\rm in},\theta_{\rm in}),
\end{aligned}\right.
\end{equation}
where the initial data $(u_{\rm in},\theta_{\rm in})$ are periodic in $x,z$
and satisfy $\nabla\cdot u_{\rm in}=0$.

The random shear displaces fluid layers at different heights relative
to one another in the spanwise direction, creating finer vertical
scales on which viscosity acts. It therefore supplies a second mixing
mechanism for the streamwise-independent modes with nonzero spanwise
frequency, which the Couette shear does not mix. Modes with nonzero
streamwise frequency are mixed by the Couette shear, while the double
means remain unmixed.

Our main result, Theorem~\ref{thm:full3d-stability}, gives sufficient
smallness conditions that admit larger perturbations of the full nonlinear
system even when the noise amplitude vanishes with viscosity. For fixed
$\beta>1/2$ and suitable $\rho$ tending to zero with $\nu$, each fixed
deterministic Sobolev perturbation within this range yields a strong
solution that is global and satisfies stability estimates on an event
of prescribed high probability. The gain follows from controlling the
streamwise-independent field generated by interactions of the nonzero
streamwise (fast) modes together with its accumulated feedback on those
modes. Our estimates retain the physical damping of the fast vertical
velocity and temperature, require no spatial localization and allow
nonzero initial double means. The precise statement is given in
Section~\ref{sec:main-result}.

\subsection{Related literature}\leavevmode\par

 \texttt{(I) Mixing and stratified Couette flow.}
The stability gain considered here builds on the connection between
mixing and enhanced dissipation. Constantin et al.~\cite{CKRZ2008} gave
a spectral criterion for relaxation enhancement, while Bedrossian and
Coti Zelati~\cite{BCZ2017} and Coti Zelati et al.~\cite{CZDE2020}
established quantitative relations between mixing and faster diffusion.
For the Navier--Stokes equations, Bedrossian et al.~\cite{BMV2016}
proved enhanced dissipation and inviscid damping in two dimensions.
The role of these decay mechanisms in nonlinear Sobolev stability
thresholds for three-dimensional Couette flow was studied by
Bedrossian et al.~\cite{BGM2017} and Wei and Zhang~\cite{WeiZhang2021}.
For stratified Couette flow, shear mixing interacts with buoyancy waves.
The classical criteria of Miles~\cite{Miles1961} and
Howard~\cite{Howard1961} concern inviscid spectral stability.
In two dimensions, Bianchini et al.~\cite{BCD2022} proved linear damping,
and Bedrossian et al.~\cite{BBCZD2023} established nonlinear inviscid
damping and shear-buoyancy growth. The viscous problem was studied by
Zillinger~\cite{Zillinger2021} and Zhai and Zhao~\cite{ZhaiZhao2023},
while Masmoudi et al.~\cite{MSHZ2022} treated the case without thermal
diffusivity. In the aligned three-dimensional geometry, Coti Zelati and
Del Zotto~\cite{CZDz2025} showed that buoyancy suppresses lift-up.
Related stability questions for rotating Couette flow were studied by
Coti Zelati et al.~\cite{CZDZWRotating2025}, Huang et
al.~\cite{HLSX2026} and Li et al.~\cite{LSWWZ2025}.

 \texttt{(II) Deterministic thresholds and the zero modes.}
For the full three-dimensional stratified system, decay of the fast
modes must be combined with control of the streamwise-independent field
and its feedback. The closest deterministic comparison is the work of
Coti Zelati et al.~\cite{CZDZW} on $\T_x\times\R_y\times\T_z$.
For fixed $\beta>1/2$, they prove full stability at size $\nu^{11/12}$
for initial data in $H^{2m+1}\cap W^{2m+5,1}$, $m\ge4$, with vanishing
initial double means of $u^3,\theta$. Their estimates retain the physical
damping powers $t^{-3/2}$ for $u^2_{\mathrm{nz}}$ and $t^{-1/2}$ for
$\theta_{\mathrm{nz}}$, together with dissipation on the $\nu^{-1/3}$
time scale. They also discuss relaxation of the mean restrictions
in~\cite{CZDZW}*{Remark~1.1}. Although the nonzero-mode estimates allow
size $\nu^{5/6}$, the zero modes and their coupling impose a smaller
radius for the full system. This restriction motivates additional mixing
of the streamwise-independent modes with nonzero spanwise frequency.

The geometry of the domain provides another way to control this field.
On $\T_x\times\R_y\times\R_z$, Chen and Zhang~\cite{ChenZhang2026}
use the unbounded spanwise direction in their dispersive estimates.
For any fixed $\beta>0$ and comparable viscosity and thermal diffusivity,
they obtain the threshold $\nu^{2/3}|\log\nu|^{-2-2a_0}$, with fixed
sufficiently small $a_0>0$. Their assumptions combine Sobolev regularity
with $W^{s,1}$ localization of the streamwise average, $s\ge4$.
At this threshold, their estimates allow weaker damping and
$\nu^{-1/6}$ high-norm growth of $u^1_0$. Under $\beta>1/2$, stronger
initial regularity and size $\nu^{3/4+\delta}$ recover the linear damping
rates for sufficiently small $\delta>0$. Since our spanwise direction
is periodic, this dispersive mechanism is unavailable. We use random
spanwise transport to obtain additional mixing and retain the linear
damping powers for the fast vertical velocity and temperature with
Sobolev data, without spatial localization.

 \texttt{(III) Transport noise and stochastic stability.}
Random transport provides a mechanism for this additional dissipation.
For prescribed stochastic shear, Coti Zelati et al.~\cite{CZHV2025}
quantified the enhancement of dissipation and established a stochastic
RAGE theorem. When the shear profiles have no common critical points,
these estimates give the $\rho\sqrt{\nu}$ rate for fixed nonzero Fourier
modes. Bedrossian et al.~\cite{BBPS2021} proved almost-sure scalar mixing
and enhanced dissipation for transport driven by stochastic
Navier--Stokes velocities. Scaling limits of transport noise were
studied by Galeati~\cite{Galeati2020} and Luo~\cite{Luo2021}, while its
regularizing effects were investigated by Flandoli and
Luo~\cite{FlandoliLuo2021}, Luo~\cite{Luo2025} and
Agresti~\cites{Agresti2025,Agresti2026}.

Stochastic stability near Couette flow has also been studied by Arbon
and Bedrossian~\cite{ArbonBedrossian2026} in two dimensions, using an
averaging principle under additive forcing and Navier boundary
conditions. For the common transport noise in the present system,
control of nonlinear sources and modal feedback must be accompanied by
a solution theory compatible with the projected transport operator.
Agresti and Veraar~\cite{AgrestiVeraar2024} developed critical-space
stochastic Navier--Stokes theory. On $\mathbb T^3$, Lin et
al.~\cite{LinLiuWang2024} proved high-probability global existence for the
Boussinesq equations with small rough data and small noise, using spatial
$L^p$ data with $p>5$ and associated probability moments.
Goodair~\cite{Goodair2025} derived closed Sobolev estimates for pure
Leray-projected transport noise on the torus. 
In our setting, the affine
coefficient $y$ is unbounded on $\mathbb R_y$. We treat this coefficient
through a projected-group construction without spatial localization,
so that the resulting solution theory can be used in the coupled
stability estimates.


\subsection{Main result}\label{sec:main-result}

\begingroup
\setlength{\parskip}{6pt plus 1pt minus 1pt}
\setlength{\abovedisplayskip}{9pt plus 2pt minus 1pt}
\setlength{\belowdisplayskip}{9pt plus 2pt minus 1pt}
\setlength{\abovedisplayshortskip}{6pt plus 2pt minus 1pt}
\setlength{\belowdisplayshortskip}{7pt plus 2pt minus 1pt}

To state the main result, let us separate the modes according to the
action of the two shears. Periodic averages are taken with respect to
normalized Haar measure.
We write $f_0$ for the $x$ average and
$f_{\mathrm{nz}}=f-f_0$ for the nonzero streamwise modes. We also write
$f_{00}$ for the $x,z$ average and set $\widetilde f_0=f_0-f_{00}$.
In Fourier variables,
$f_{\mathrm{nz}}$, $\widetilde f_0$, and $f_{00}$ correspond to
$k\ne0$, $(k=0,\ell\ne0)$, and $(k=\ell=0)$, respectively.
The divergence constraint and $L^2_y$ integrability imply $u^2_{00}=0$.

We introduce the volume-preserving map
\begin{equation}\label{eq:pullback}
 \begin{aligned}
 \Phi_t(X,Y,Z)&=(X+tY,Y,Z+\rho YW_t),\\[2pt]
 D_t&=(\partial_X,\partial_Y-t\partial_X-\rho W_t\partial_Z,\partial_Z).
 \end{aligned}
\end{equation}
Set
\[
 f^\sh=f\circ\Phi_t, \quad \nabla_t=D_t, \quad
 \|f\|_{H^s_{\ad(t)}}=\|f^\sh\|_{H^s}.
\]
In physical Fourier variables this norm has weight
$\av{k,\eta+kt+\rho\ell W_t,\ell}^{2s}$, where
$\av{\xi}=(1+|\xi|^2)^{1/2}$.
The velocity is pulled back componentwise. The projected group
constructed in Section~\ref{sec:local} also rotates the divergence-free
Fourier plane and preserves the transformed Sobolev norms. We use the
strong solution class of Proposition~\ref{thm:affine-local-strong},
defined through this group or, equivalently, through the weak
Stratonovich formulation against smooth compactly supported test functions.

For a Sobolev index $s$, set
\[
 v_0=(u_0^2,u_0^3), \quad q_0=u_0^1+\beta^{-1}\theta_0, \quad
 Y_0=(v_0,\theta_0,q_0).
\]
We define
\[
 \begin{aligned}
 E_0=\|Y_0^\sh\|_{H^s}^2,
 \quad G_0=\|\nabla_tY_0^\sh\|_{H^s}^2,\quad
 S=\|\widetilde{Y_0^\sh}\|_{H^s},
 \quad G_{\rm s}=\|\nabla_t\widetilde{Y_0^\sh}\|_{H^s}^2.
 \end{aligned}
\]

An overbar on $Y_0$ denotes its $z$-average, and a tilde signifies the remainder obtained after subtracting this average. Thus $E_0$ incorporates both zero sectors.
We omit the identically zero vertical velocity component from mean vectors. We separate the initial mean from the nonlinear response by setting
\[
m_{\rm in}=\overline{Y_0}(0), \quad
m_{\rm nl}(t)=\overline{Y_0}(t)-e^{\nu t\partial_{yy}}m_{\rm in}.
\]

For the fast modes, we use the reference Fourier variables
$\xi=(k,\eta,\ell)$. For $k\ne0$, put
\[
 \begin{aligned}
 h=(k^2+\ell^2)^{1/2}, \quad \eta_t=\eta-kt-\rho\ell W_t,\quad
 r_t=(h^2+\eta_t^2)^{1/2}, \quad \lambda_k=\nu^{1/3}|k|^{2/3}.
 \end{aligned}
\]
The divergence-free plane has the orthonormal basis
\begin{equation}\label{eq:intro-polarization}
 e_g(\eta_t)=\frac{(-k\eta_t,h^2,-\ell\eta_t)}{hr_t},
  \quad e_v=\frac{(\ell,0,-k)}h.
\end{equation}
Write
\[
 \begin{aligned}
 \widehat{u^\sh} =e_g(\eta_t)v+e_vb, \quad
 p=b+\frac{\ell}{\beta h}\widehat{\theta^\sh},\quad
 J=r_t\bigl(|v|^2+|\widehat{\theta^\sh}|^2\bigr)+|p|^2.
 \end{aligned}
\]
For the fixed Sobolev index $s$, define the fast-mode energies
\begin{equation}\label{eq:intro-fast-energies}
 \begin{aligned}
  E_{\mathrm{nz}}(t)&=\sum_{k\ne0,\ell}\int_\R
                         \av{\xi}^{2s}J(t,\xi)\dd\eta,\quad
  D_{\mathrm{nz}}(t) =\sum_{k\ne0,\ell}\int_\R
                         \av{\xi}^{2s}r_t^2J(t,\xi)\dd\eta,\\[2pt]
  L_{\mathrm{nz}}(t)&=\sum_{k\ne0,\ell}\int_\R
                         \av{\xi}^{2s}\lambda_kJ(t,\xi)\dd\eta.
 \end{aligned}
\end{equation}
For fixed $\beta>1/2$, these quantities are uniformly equivalent
to the compensated energies in \eqref{eq:lp-fast-field-energy}.
They control the adapted $H^s$ norm of the fast field, and initial
$H^{s+1/2}$ regularity makes $ E_{\mathrm{nz}}(0)$ finite. The weights
retain the different derivative orders of the wave and passive components.

We are now in a position to state the main result of this paper.

\endgroup

\Needspace{8\baselineskip}
\begingroup
\setlength{\parskip}{5pt plus 1pt minus 1pt}
\setlength{\jot}{7pt}
\begin{theorem}\label{thm:full3d-stability}
Fix an integer $s\ge3$ and $\beta>1/2$. There is a constant
$c_{s,\beta}>0$ such that the following holds. Let $0<\nu\le1$ and
$\rho\ge0$ satisfy
\[
 \rho^2\le c_{s,\beta}\nu^{1/3}, \quad
 \nu\le\mu=\nu+\rho\sqrt\nu\le\lambda=\nu^{1/3}.
\]
Let $(u_{\rm in},\theta_{\rm in})\in H^s(\Omega;\mathbb R^3)\times H^s(\Omega)$
be deterministic, with $\nabla\cdot u_{\rm in}=0$ and
$E_{\mathrm{nz}}(0)<\infty$. Set
\[
 \begin{aligned}
 \eps_0^2=E_0(0), \quad \alpha^2=E_{\mathrm{nz}}(0),\quad
 \eps^2=\eps_0^2+\alpha^2,
 \quad \delta=\eps_0+\frac{\alpha^2}{\sqrt{\nu\lambda}}.
 \end{aligned}
\]
For every $\zeta\in(0,1)$, if
\begin{equation}\label{eq:full3d-small}
 \begin{aligned}
 \alpha \le c_{s,\beta}\zeta^{1/4}
        \min\{\nu^{3/4},\nu^{7/12}\mu^{1/4}\},\quad
 \eps_0 \le c_{s,\beta}\zeta^{1/4}\,
                  \min\{\nu^{5/6},\sqrt{\nu\mu}\},
 \end{aligned}
\end{equation}
then \eqref{eq:model}, driven by $W$ with initial data
$(u_{\rm in},\theta_{\rm in})$, admits a unique maximal strong solution
$(u,\theta)$ in the class of Proposition~\ref{thm:affine-local-strong}.
Writing $\tau_*$ for its maximal lifetime, there is an event
$\Gamma_\zeta$ with $\Prob(\Gamma_\zeta)\ge1-\zeta$ such that
$\tau_*=\infty$ on $\Gamma_\zeta$. On this event, the solution is global
and all the following conclusions hold. The constants
$C_{s,\beta,\zeta}$ below depend only on $s,\beta,\zeta$, and the Sobolev
energy estimates use the reference norms defined above.

\begin{enumerate}[label=(\roman*),leftmargin=*,itemsep=8pt,topsep=6pt]
 
\item \textup{\texttt{Energy bounds and enhanced dissipation.}}
The zero and fast fields satisfy
\begin{align}
 \sup_{t\ge0}(E_0+E_{\mathrm{nz}})
       +\int_0^\infty[\nu(G_0+D_{\mathrm{nz}})+L_{\mathrm{nz}}]\dd t
 &\le C_{s,\beta,\zeta}\eps^2,\label{eq:full3d-total}\\
 \sup_{t\ge0}E_0(t)+\nu\int_0^\infty G_0\dd t
 &\le C_{s,\beta,\zeta}\delta^2.\label{eq:full3d-zero-budget}
\end{align}
For some $b=b_{s,\beta}>0$, the enhanced-dissipation estimates are
\begin{align}
 \sup_{t\ge0}e^{\lambda t/8}E_{\mathrm{nz}}(t)
       +\int_0^\infty e^{\lambda t/8}(\nu D_{\mathrm{nz}}+L_{\mathrm{nz}})\dd t
 &\le C_{s,\beta,\zeta}\alpha^2,\label{eq:full3d-fast}\\
 \sup_{t\ge0}e^{b\mu t}S(t)^2
       +\nu\int_0^\infty e^{b\mu t}G_{\rm s}(t)\dd t
 &\le C_{s,\beta,\zeta}\delta^2.\label{eq:full3d-simple}
\end{align}
The nonlinear response of the double mean satisfies
\begin{equation}\label{eq:full3d-mean}
 \sup_{t\ge0}\norm{m_{\rm nl}(t)}_{H^s}^2
       +\nu\int_0^\infty\norm{\partial_y m_{\rm nl}}_{H^s}^2\dd t
 \le C_{s,\beta,\zeta}
       \left(\frac{\delta^4}{\nu\mu}
                         +\frac{\alpha^4}{\nu\lambda}\right).
\end{equation}
The same estimate holds with $m_{\rm nl}$ replaced by
$\overline{Y_0}$ and $\|m_{\rm in}\|_{H^s}^2$ added to the right-hand
side, after changing $C_{s,\beta,\zeta}$.
 
\item \textup{\texttt{Physical  damping and strong convergence.}}
For every $t\ge0$, the physical fast components satisfy
\begin{equation}\label{eq:full3d-physical-damping}
 \begin{aligned}
 \norm{u_{\mathrm{nz}}^2(t)}_{L^2(\Omega)}
 &\le C_{s,\beta,\zeta}\alpha\av{t}^{-3/2}e^{-\lambda t/32},\\[2pt]
 \norm{\theta_{\mathrm{nz}}(t)}_{L^2(\Omega)}
 &\le C_{s,\beta,\zeta}\alpha\av{t}^{-1/2}e^{-\lambda t/32}.
 \end{aligned}
\end{equation}
As $t\to\infty$, $\overline{Y_0}(t)$ converges strongly to zero in
$H^s(\R)$, and both $(u^\sh,\theta^\sh)(t)$ and $(u,\theta)(t)$ converge
strongly to zero in $H^s(\Omega)\times H^s(\Omega)$.

\item \textup{\texttt{Admissible perturbation scales.}}
The common-size condition
\[
 \eps\le c_{s,\beta}\zeta^{1/4}\min\{\nu^{5/6},\sqrt{\nu\mu}\}
\]
is sufficient. For fixed $\zeta$ and the noise amplitudes below, the
admissible fast- and zero-mode orders are
\begin{center}
\begin{tabular}{@{}ccc@{}}
\toprule
\textup{Noise amplitude} & \textup{Order of $\alpha$} & \textup{Order of $\eps_0$}\\
\midrule
$\rho=a\nu^{1/6}$ & $\nu^{3/4}$ & $\nu^{5/6}$\\
$\rho=\nu^{1/4}$ & $\nu^{37/48}$ & $\nu^{7/8}$\\
$\rho=0$ & $\nu^{5/6}$ & $\nu$\\
\bottomrule
\end{tabular}
\end{center}
Here $a>0$ is fixed and sufficiently small, and the first row applies
for sufficiently small $\nu$. The fast scales also apply to initially
streamwise-mean-free data ($\eps_0=0$), with the generated zero field
included in $\delta$.
\end{enumerate}
\end{theorem}
\endgroup

\begin{remark}
For fixed $\zeta\in(0,1)$ and $\rho=a\nu^{1/6}$, with fixed sufficiently
small $a>0$ and sufficiently small $\nu$, Theorem~\ref{thm:full3d-stability}
admits zero- and fast-mode
sizes of order $\nu^{5/6}$ and $\nu^{3/4}$ while retaining the physical
damping powers $t^{-3/2}$ and $t^{-1/2}$. For data satisfying the
assumptions of both results, our sufficient common-size radius $\nu^{5/6}$
exceeds the fixed-$\beta$ sufficient radius $\nu^{11/12}$ established
in~\cite{CZDZW} by a factor $\nu^{-1/12}$, even though $\rho\to0$.
The event of probability at least $1-\zeta$ may depend on the
deterministic initial datum.
\end{remark}
\begin{remark}
For integers $s\ge3$, Theorem~\ref{thm:full3d-stability}
applies to data with wave amplitudes in $H^{s+1/2}$ and passive and zero
components in $H^s$, without additional spatial localization or $W^{k,1}$
assumptions. Nonzero initial double means are also allowed within the
zero-sector smallness bound.
\end{remark}
\Needspace{8\baselineskip}
\begin{remark}
The novelty lies in coupling noise-enhanced simple-zero decay with control
of fast self-interactions and accumulated feedback.
Lemma~\ref{lem:compensated-fast} applies to other coupled shear-flow systems
under its hypotheses. Extending the full method requires compatible modal
energies, a bounded frequency correction, conditional divergence smoothing
and nonlinear estimates that also control the unmixed means.
\end{remark}
 
\Needspace{10\baselineskip}
The temperature double mean admits a uniform quadratic estimate relative
to its free heat evolution. In the following corollary, an overbar
 denotes normalized $x,z$ averaging.

\begin{corollary}\label{cor:full-temperature-mean}
Assume the hypotheses of Theorem~\ref{thm:full3d-stability} at
index $s+2$, with $\eps$ defined by the joint initial energy at
that index. Then, on the event of that theorem,
\begin{equation}\label{eq:full-temperature-mean}
 \sup_{t\ge0}\|\overline\theta(t)-e^{\nu t\partial_{yy}}\overline{\theta_{\rm in}}\|_{H^s(\R)}
          \le C_{s,\beta,\zeta}\eps^2.
\end{equation}
The constant is uniform in $\nu,\rho$ in the stated range.
For example, $H^{s+5/2}$ initial regularity suffices.
\end{corollary}

\subsection{Difficulties and main ideas of the proof}\label{rem:proof-ideas}

The main difficulty is to pass from the additional mixing supplied
by the random spanwise shear to stability of the full nonlinear system.
Fast-mode interactions generate a zero field, and this field feeds back
into the fast equations. The random shear mixes the simple-zero modes,
but the double means have no spectral gap. Two further features enter
the estimates: the nonlinear fluxes are adapted to the same noise, and
the feedback contains a physical half derivative. The argument must
control these interactions without losing the damping of the fast
vertical velocity and temperature.

To this end, we use noise-enhanced simple-zero decay to control both the
generated zero field and its feedback on the fast modes. This is the
mechanism that yields a larger stability range for suitable noise
amplitudes vanishing with viscosity, while allowing nonzero initial
double means. To implement it, we combine a bounded Poisson correction
with conditional estimates for nonlinear fluxes and an adapted weight
for the accumulated feedback. The roles of these constructions are
described in the following steps.

 \texttt{Step 1 (Projected transport and local strong solutions).}
We begin with the projected transport generated by the unbounded
affine coefficient $y$. Fourier polarization gives a transport group
that preserves the divergence-free constraint, with both Leray
projections retained in the It\^o correction. By the cancellations in
Lemma~\ref{lem:affine-noise-energy}, we obtain Sobolev bounds without
spatial localization. These bounds are then used to establish the local
strong solution theory in Proposition~\ref{thm:affine-local-strong}.

 \texttt{Step 2 (Critical frequencies and fast self-interactions).}
The wave and passive components of the fast field carry different
derivative weights. We first separate their nonlinear interactions
using the Fourier structure in Lemma~\ref{lem:separated-mixed}.
A loss remains near critical frequencies. To control it, we introduce
a bounded Poisson correction for the generator of the stochastic
frequency characteristics (Lemma~\ref{lem:critical-poisson}). The
resulting positive term $Z_g$ absorbs $C\nu^{-3}E_g^2Z_g$ when the
fast size is a sufficiently small multiple of $\nu^{3/4}$. The remaining
noise drift is absorbed under $\rho^2\le c_{s,\beta}\lambda$.

 \texttt{Step 3 (The zero-field feedback on the fast modes).}
On a common stopping interval, we impose a bound on the total
zero-mode energy and dissipation. The feedback coefficient contains
a physical half derivative measured in the reference $H^s$ norm.
Energy--dissipation interpolation controls its time integral, which we
incorporate into an adapted time weight. Under
$\delta\lesssim\nu^{5/6}$, fast-mode decay absorbs the accumulated loss
(Lemma~\ref{lem:compensated-fast}). We thus obtain the fast-mode moments
required to estimate the zero field generated by fast interactions.

 \texttt{Step 4 (Conditional response of the simple-zero modes).}
The flux at time $r$ depends on the noise up to $r$, so we estimate
its response conditionally. We apply the divergence smoothing
estimate for the coupled projected system in Lemma~\ref{lem-smooth}.
The kernel $C_p[\nu(t-r)]^{-1/2}e^{-c_p\mu(t-r)}$ controls the physical
derivative at the forcing time and has time integral of size
$(\nu\mu)^{-1/2}$. The simple-zero nonlinear feedback can therefore be
absorbed under $\delta\lesssim\sqrt{\nu\mu}$.
The fast-mode moments from Step 3, together with
Lemma~\ref{lem:fast-source-convolution}, bound the response generated
by the fast modes by $\alpha^2/\sqrt{\nu\lambda}$. This contribution
is included in the total zero-mode budget, so the smallness condition
accounts for the generated field as well as the initial perturbation.

\texttt{Step 5 (The double means).}
The double means have no spectral gap and require a separate
estimate. We split off their free heat evolution, retaining the
contribution of the initial means and estimating the response generated
by the nonlinear fluxes.

 \texttt{Step 6 (Closing the coupled energy estimates).}
We now close the fast and total zero-mode estimates using the
common stopping time, on an event of probability at least $1-\zeta$.
For prescribed failure probability $\zeta$, the three restrictions take
the form
\[
 \begin{aligned}
  \alpha&\lesssim\zeta^{1/4}\nu^{3/4}
       &&\text{(fast self-interactions)},\\
  \delta&\lesssim\zeta^{1/4}\nu^{5/6}
       &&\text{(accumulated zero-field feedback)},\\
  \delta&\lesssim\zeta^{1/4}\sqrt{\nu\mu}
       &&\text{(simple-zero kernel absorption)},
 \end{aligned}
\]
where the implicit constants depend only on $s,\beta$. With the generated
zero field included in $\delta$, these restrictions lead to the sufficient
conditions in~\eqref{eq:full3d-small} and preserve fast-mode decay at rate
$\lambda=\nu^{1/3}$. Finally, the energy bounds and control of the
Brownian frequency shift give the physical damping of the fast vertical
velocity and temperature. This completes the global stability argument.

\subsection{Organization of the paper}

Sections~\ref{sec:local} and~\ref{sec:lp-linear} provide the local
strong solution theory and conditional linear estimates, respectively.
Using these estimates, Section~\ref{sec:nonlinear} treats the
streamwise-independent dynamics and establishes nonlinear decay and
scalar mean normal forms. Section~\ref{sec:3d-stability} combines the
fast- and zero-mode estimates to prove Theorem~\ref{thm:full3d-stability}
and Corollary~\ref{cor:full-temperature-mean}.

\setlength{\parskip}{5pt plus 1pt minus 1pt}
\setlength{\jot}{7pt}
\setlength{\abovedisplayskip}{11pt plus 2pt minus 2pt}
\setlength{\belowdisplayskip}{11pt plus 2pt minus 2pt}
\setlength{\abovedisplayshortskip}{7pt plus 2pt minus 1pt}
\setlength{\belowdisplayshortskip}{9pt plus 2pt minus 1pt}
\makeatletter
\def\thm@space@setup{%
  \thm@preskip=10pt plus 2pt minus 1pt
  \thm@postskip=10pt plus 2pt minus 1pt}
\makeatother
\AddToHook{env/proof/before}{\Needspace{4\baselineskip}}

\section{Projected affine transport and local solutions}\label{sec:local}

The affine coefficient is unbounded in the vertical direction. We first
derive the cancellation that permits Sobolev energy estimates without
spatial weights, and then use the projected transport group to construct
local strong solutions.
 Denote by $H^s_\sigma(\Omega)$ the divergence-free subspace of
$H^s(\Omega;\mathbb R^3)$ and by $\mathcal P$ the Leray projection. For Fourier
variables $(k,\eta,\ell)\in\mathbb Z\times\mathbb R\times\mathbb Z$, set
$A=\mathcal P(y\partial_x)$, $N=\mathcal P(y\partial_z)$,
$A_\theta=y\partial_x$, and $N_\theta=y\partial_z$.
We use $\Lambda=(1-\Delta)^{1/2}$ and restrict velocity operators to
$H^s_\sigma$. In particular, $N^2=\mathcal P(y\partial_z)\mathcal P(y\partial_z)$ differs
from $\mathcal P(y^2\partial_z^2)$.

\begin{lemma}\label{lem:affine-noise-energy}
Let $s\geq0$. For every smooth divergence-free $u$ in a common operator core,
\begin{align*}
 \big|\langle u,N^2u\rangle_{H^s}+\|Nu\|_{H^s}^2\big|
 &\leq K_s\|u\|_{H^s}^2,\\
 |\langle u,Nu\rangle_{H^s}|&\leq\frac{s}{2}\|u\|_{H^s}^2,\end{align*}
where $K_s=s+\tfrac12s|s-1|$. The same bounds hold for $N_\theta$.
After cancellation, the expressions on the left extend continuously to
$H^s_\sigma$ and $H^s$, respectively. The constants are independent of
the range of $y$.
\end{lemma}

\begin{proof}
For $h=(k^2+\ell^2)^{1/2}>0$ and $r=(h^2+\eta^2)^{1/2}$, define
\[
 e_g=\frac{(-k\eta,h^2,-\ell\eta)}{hr}, \quad
 e_v=\frac{(\ell,0,-k)}{h}, \quad E=(e_g,e_v), \quad
 n=\frac{(k,\eta,\ell)}r.
\]
These orthonormal vectors satisfy $E^*E=I_2$, $EE^*=\mathcal P$,
$\partial_\eta e_g=-\frac{h}{r^2}n$, $\partial_\eta e_v=0$, and
$\mathcal P\partial_\eta E=0$.
Writing $\widehat u=Ea$, $a:\mathbb R\to\mathbb C^2$, and using the fact that $y\partial_z$
corresponds to $-\ell\partial_\eta$, we obtain after projection
$\widehat{Nu}=-\ell E\partial_\eta a$ and
$\widehat{N^2u}=\ell^2E\partial_\eta^2a$, including the projection
commutators. On the modes $h=0$, both $A$ and $N$ vanish.
Set $w_s=(1+k^2+\eta^2+\ell^2)^s$.
Integration by parts yields
\begin{align}
 \langle u,N^2u\rangle_{H^s}+\|Nu\|_{H^s}^2
 &=\frac12\sum_{h>0}\int_{\mathbb R}\ell^2w_s''|a|^2\,d\eta,
 \label{eq:affine-drift-exact}\\
 \langle u,Nu\rangle_{H^s}
 &=\frac12\sum_{h>0}\int_{\mathbb R}\ell w_s'|a|^2\,d\eta.
 \label{eq:affine-mart-exact}
\end{align}
All vector pairings denote real parts. Direct differentiation gives
\[
 \frac{w_s'}{w_s}=\frac{2s\eta}{1+r^2}, \quad
 \frac{w_s''}{w_s}=\frac{2s}{1+r^2}
       +\frac{4s(s-1)\eta^2}{(1+r^2)^2}.
\]
Using $2|\ell\eta|\leq1+r^2$ and
$4\ell^2\eta^2\leq(1+r^2)^2$, we obtain
$\frac{|\ell w_s'|}{w_s}\leq s$ and
$\frac{|\ell^2w_s''|}{w_s}\leq2s+s|s-1|$.
The scalar estimates follow from the same calculation, with $E$ omitted.
\end{proof}

We can also express this cancellation in commutator form. On $H^s_\sigma$, set
$C_s=[\Lambda^s,N]\Lambda^{-s}$. Since $\Lambda^s\mathcal P=\mathcal P\Lambda^s$,
$C_s$ is the real scalar Fourier multiplier
$c_s(k,\eta,\ell)=\frac{s\ell\eta}{1+k^2+\eta^2+\ell^2}$.
If $v=\Lambda^su$, then $\Lambda^sNu=(N+C_s)v$ and
$\Lambda^sN^2u=(N+C_s)^2v$.  As $N^*=-N$ and $C_s^*=C_s$,
\begin{equation*}
 \langle v,(N+C_s)^2v\rangle+\|(N+C_s)v\|_2^2
 =\langle v,([C_s,N]+2C_s^2)v\rangle.
\end{equation*}
The right-hand side has symbol
\[
 \ell\partial_\eta c_s+2c_s^2
 =\frac{s\ell^2}{1+r^2}
      +\frac{2s(s-1)\ell^2\eta^2}{(1+r^2)^2}
 =\frac{\ell^2w_s''}{2w_s}.
\]
Thus, even after taking high-order commutators, the It\^o correction
and quadratic variation together contribute only a bounded zeroth-order
term. Applied to \eqref{eq:model}, this cancellation gives the energy
estimate below.

\begin{proposition}\label{prop:affine-local-energy}
Let $s>5/2$. For a sufficiently regular solution of
\eqref{eq:model}, set
\[
 \mathcal E_s^{\mathrm{phys}}=\|u\|_{H^s}^2+\|\theta\|_{H^s}^2,
  \quad
 \mathcal D_s^{\mathrm{phys}}=\|\nabla u\|_{H^s}^2+\|\nabla\theta\|_{H^s}^2.
\]
Then we have
\begin{equation}\label{eq:affine-total-energy}
 \begin{aligned}
 d\mathcal E_s^{\mathrm{phys}}+2\nu\mathcal D_s^{\mathrm{phys}}\,dt
 &\leq C_s\bigl(1+\rho^2+\|\nabla u\|_\infty
                        +\|\nabla\theta\|_\infty\bigr)
                   \mathcal E_s^{\mathrm{phys}}\,dt+m_s\,dW_t,\\
 |m_s|&\leq s|\rho|\mathcal E_s^{\mathrm{phys}}.
 \end{aligned}
\end{equation}
The coefficient of $\rho^2\mathcal E_s^{\mathrm{phys}}$ depends only on $s$,
uniformly in $\nu$, $\beta$, spatial cutoffs, and $y$.
\end{proposition}

\begin{proof}
The velocity It\^o correction and quadratic variation give
$\rho^2\bigl(\langle u,N^2u\rangle_{H^s}+\|Nu\|_{H^s}^2\bigr)dt$.
Adding the temperature contribution, we obtain the bound
$\rho^2K_s\mathcal E_s^{\mathrm{phys}}dt$. The martingale coefficient
\[
 m_s=-2\rho\bigl(\langle u,Nu\rangle_{H^s}
       +\langle\theta,N_\theta\theta\rangle_{H^s}\bigr)
\]
is controlled by Lemma~\ref{lem:affine-noise-energy}. Replacing $\ell$
by $k$ in \eqref{eq:affine-mart-exact} bounds the affine drift
by $s\mathcal E_s^{\mathrm{phys}}$. The stretching term is bounded by
$\mathcal E_s^{\mathrm{phys}}$, and the buoyancy terms cancel:
$-2\beta\langle u,\mathcal P(\theta e_y)\rangle_{H^s}
+2\beta\langle\theta,u^2\rangle_{H^s}=0$.
For the nonlinear terms, incompressibility and the Sobolev commutator estimate give
\begin{align*}
 |\langle\Lambda^su,\Lambda^s((u\cdot\nabla)u)\rangle|
 &\leq C_s\|\nabla u\|_\infty\|u\|_{H^s}^2,\\
 |\langle\Lambda^s\theta,\Lambda^s((u\cdot\nabla)\theta)\rangle|
 &\leq C_s\bigl(\|\nabla u\|_\infty\|\theta\|_{H^s}^2
       +\|\nabla\theta\|_\infty\|u\|_{H^s}\|\theta\|_{H^s}\bigr).
\end{align*}
Combining these bounds, we obtain \eqref{eq:affine-total-energy}.
\end{proof}

\begin{proposition}\label{thm:affine-local-strong}
Let $s>5/2$, $\nu>0$, and $\rho,\beta\in\mathbb R$. On a filtered probability space carrying a real Brownian motion,
assume that $(u_{\rm in},\theta_{\rm in})$ is $\mathcal F_0$-measurable and almost surely in $H^s_\sigma\times H^s$. Then \eqref{eq:model} admits a unique maximal adapted solution with an almost surely positive lifetime $\tau_*$ such that, for every stopping time $\tau<\tau_*$,
\[
 (u,\theta)\in C([0,\tau];H^s_\sigma\times H^s)
       \cap L^2(0,\tau;H^{s+1}_\sigma\times H^{s+1})
 \quad\hbox{almost surely}.
\]
The solution is probabilistically strong and pathwise unique.
If $\tau_*<\infty$, then
\[
 \limsup_{t\uparrow\tau_*}
 \bigl(\|u(t)\|_{H^s}+\|\theta(t)\|_{H^s}\bigr)=\infty.
\]
Moreover, finite-time breakdown requires
\[
 \int_0^{\tau_*}
 \bigl(\|\nabla u(t)\|_\infty+\|\nabla\theta(t)\|_\infty\bigr)dt=\infty.
\]
We understand the equation through the transport group constructed below
or, equivalently, in distributional It\^o form. For general unweighted
$H^s$ initial data, this formulation does not require an $H^{s-1}$-valued It\^o integral.
\end{proposition}

\Needspace{4\baselineskip}
\begin{proof}
For $h>0$, the amplitude $a=E^*\widehat u$ identifies the divergence-free
Fourier fibre with $\mathbb C^2$. In these coordinates, $A=-k\partial_\eta$ and
$N=-\ell\partial_\eta$ generate commuting unitary $L^2$ groups:
\[
 \widehat{R(a_0,b_0)u}(k,\eta,\ell)
   =E(k,\eta,\ell)a(k,\eta+ka_0+\ell b_0,\ell),
  \quad R(a_0,b_0)=e^{-a_0A-b_0N}.
\]
For $h=0$, set $R=I$. The scalar group uses the same Fourier shift, with no frame rotation. For every real $r$ and finite $a_0,b_0$,
\begin{equation}\label{eq:affine-group-hs}
 \|R(a_0,b_0)^{\pm1}\|_{H^r\to H^r}
       \leq C_r(1+|a_0|+|b_0|)^{|r|}.
\end{equation}
Comparison of the Fourier weights proves this bound; the groups are
also strongly continuous on $H^r$. We now remove the affine transport by setting
$\mathcal R_t=\operatorname{diag}(R(t,\rho W_t),R_\theta(t,\rho W_t))$
and $V=\mathcal R_t^{-1}(u,\theta)$.
The Stratonovich chain rule cancels both transport terms. The conjugated
Laplacian has symbol $-q_t(k,\eta,\ell)$, where
$q_t=k^2+(\eta-kt-\rho\ell W_t)^2+\ell^2$.
For each finite $T$, continuity of $W$ gives random constants
$0<c_T\le C_T<\infty$ such that
\[
 c_T(k^2+\eta^2+\ell^2)\leq q_t
        \leq C_T(k^2+\eta^2+\ell^2), \quad 0\leq t\leq T.
\]
The diffusion evolution family
$\mathcal S(t,r)=\exp(-\nu\int_r^tq_a\,da)$, $t\ge r$, is bounded on
$H^s$ and satisfies, for $0<t-r\leq T$,
\[
 \|\mathcal S(t,r)F\|_{H^s}
 \leq C_T\bigl[1+(\nu(t-r))^{-1/2}\bigr]\|F\|_{H^{s-1}}.
\]

Let $\mathcal L$ denote the stretching and buoyancy terms, and $\mathcal Q$ the
quadratic convection in \eqref{eq:model}. Setting
$\mathcal L_t=\mathcal R_t^{-1}\mathcal L\mathcal R_t$ and
$\mathcal Q_t=\mathcal R_t^{-1}\mathcal Q(\mathcal R_t\cdot,
\mathcal R_t\cdot)$, we have
\[
 \partial_tV=\nu\mathcal R_t^{-1}\Delta\mathcal R_tV
        -\mathcal L_tV-\mathcal Q_t(V,V).
\]
On each finite interval, \eqref{eq:affine-group-hs} and the Sobolev product estimate give
\begin{align*}
 \|\mathcal L_tV\|_{H^s}&\leq C_{s,T,\beta}\|V\|_{H^s},\\
 \|\mathcal Q_t(V,V)\|_{H^{s-1}}&\leq C_{s,T}\|V\|_{H^s}^2,\\
 \|\mathcal Q_t(V,V)-\mathcal Q_t(Z,Z)\|_{H^{s-1}}
 &\leq C_{s,T}(\|V\|_{H^s}+\|Z\|_{H^s})\|V-Z\|_{H^s}.
\end{align*}
The integrability of $(t-r)^{-1/2}$ makes the Duhamel map a contraction
on a sufficiently short interval, yielding local existence, uniqueness,
and continuous dependence. To restart at a finite stopping time $\sigma$,
we use
\[
 \mathcal R_{\sigma,t}=\operatorname{diag}
 \bigl(R(t-\sigma,\rho(W_t-W_\sigma)),
       R_\theta(t-\sigma,\rho(W_t-W_\sigma))\bigr)
\]
and stop at $\inf\{t\ge\sigma:|\rho(W_t-W_\sigma)|\ge1\}$.
On each $\mathcal F_\sigma$-measurable initial-data ball where
$\|u(\sigma)\|_{H^s}+\|\theta(\sigma)\|_{H^s}\le n$, the stopped
coefficients admit a deterministic contraction duration
$\delta_n\in(0,1]$, chosen nonincreasing in $n$.
This construction is used up to the minimum of the relative Brownian
exit time and $\sigma+\delta_n$.
The iterates at time $t$ depend only on the data at $\sigma$ and the
Brownian increments up to $t$, so the restarted solution is adapted.
If the physical $H^s$ norm stays bounded on a finite interval, the
contraction durations have a positive lower bound. The relative
Brownian exits also cannot accumulate, by uniform continuity of $W$.
Starting with $\sigma=0$ and iterating gives the adapted maximal solution
and permits continuation whenever the physical $H^s$ norm remains bounded.

We next combine the transformed diffusion energy with
\[
 |\langle V,\mathcal Q_t(V,V)\rangle_{H^s}|
 \leq C_{s,T}\|V\|_{H^{s+1}}\|V\|_{H^s}^2
 \leq\frac{\nu c_T}{4}\|\nabla V\|_{H^s}^2
         +C_{s,T,\nu}(\|V\|_{H^s}^2+\|V\|_{H^s}^4)
\]
to obtain $L^2H^{s+1}$ regularity on intervals where the solution is bounded.
The inverse group is locally bounded and transfers this regularity to
the physical variables. Restarting the construction proves the blow-up
alternative.

It remains to justify the distributional It\^o equation and the energy
estimate. We first use smooth Fourier-core approximation. For divergence-free Schwartz tests, $N\phi,N^2\phi\in L^2$;
hence the velocity noise is $\rho\langle u,N\phi\rangle dW_t$, with
correction $(\rho^2/2)\langle u,N^2\phi\rangle dt$. The scalar equation
is treated in the same way, without assuming $y\partial_zu\in H^{s-1}$.
For the energy identity, we retain the cancellation between
$\|Nu\|_{H^s}^2$ and $\langle u,N^2u\rangle_{H^s}$, which need not exist
separately. In group coordinates, the physical energy equals
\[
 \sum_{k,\ell}\int_{\mathbb R}
 w_s(k,\eta-kt-\rho\ell W_t,\ell)|a_V(t,k,\eta,\ell)|^2\,d\eta,
\]
plus its scalar analogue. After Fourier truncation, the
amplitudes have finite variation, so we apply It\^o's formula to the
weight. Its first two
Brownian derivatives are $-\rho\ell w_s'$ and $\rho^2\ell^2w_s''$, which give
\eqref{eq:affine-drift-exact}--\eqref{eq:affine-mart-exact} in physical
frequencies. On bounded time and solution intervals, the weights are
equivalent, with bounded normalized derivatives. The
$L^2H^{s+1}$ regularity controls the diffusion and convection pairings.
Passing to the limit in the combined identity gives
\eqref{eq:affine-total-energy}.

Finally, on $\mathcal E_s^{\mathrm{phys}}>0$ write
$m_s=g_s\mathcal E_s^{\mathrm{phys}}$, with $|g_s|\leq s|\rho|$.
Applying It\^o's formula to $\log(\mathcal E_s^{\mathrm{phys}}+\varepsilon)$,
discarding the nonpositive quadratic-variation term, and letting
$\varepsilon\downarrow0$ bounds finite-time logarithmic growth by
\[
 C_s(1+\rho^2)t
 +C_s\int_0^t(\|\nabla u\|_\infty+\|\nabla\theta\|_\infty)\,da
 +\int_0^tg_s\,dW_a.
\]
The martingale has bounded quadratic variation on finite intervals,
so this logarithmic bound proves the continuation criterion.
\end{proof}

\begin{lemma}\label{lem:mixed-local}
Let $s\ge3$, $\nu>0$, $\rho\ge0$, and $\beta>1/2$. Assume that
$E_0(0)+ E_{\mathrm{nz}}(0)<\infty$ almost surely. For the solution given by
Proposition~\ref{thm:affine-local-strong} and every finite stopping time
$\tau<\tau_*$, the amplitudes in \eqref{eq:intro-polarization} satisfy,
in the reference Fourier variables,
\[
\begin{gathered}
 (v,\widehat{\theta_{\mathrm{nz}}^\sh})\in
 C([0,\tau];H^{s+1/2})\cap L^2(0,\tau;H^{s+3/2}),
 \\
 (p,Y_0^\sh)\in C([0,\tau];H^s)\cap L^2(0,\tau;H^{s+1})
\end{gathered}
\]
almost surely. It follows that $ E_{\mathrm{nz}}$ is continuous and that the nonlinear
polarization energy identities hold up to finite locally bounded stopping
times.
\end{lemma}

\begin{proof}
We stop when $t$, $|\rho W_t|$, or the physical $H^s$ norm exceeds a fixed
bound. On the resulting interval, the group norms and frequency weights are equivalent.
The local solution bounds and interpolation give
$U\in C_tH^s\cap L_t^2H^{s+1}$ and $U\in L_t^4H^{s+1/2}$,
where $U=(u,\theta)$.
The transformed nonlinear forcing $\mathcal F=-\mathcal R_t^{-1}
\mathcal Q(U,U)$ therefore satisfies
\[
 \|\mathcal F(t)\|_{H^{s-1}}\le C\|U(t)\|_{H^s}^2,
  \quad
 \int_0^\tau\|\mathcal F(t)\|_{H^{s-1/2}}^2\dd t
 \le C\int_0^\tau\|U(t)\|_{H^{s+1/2}}^4\dd t<\infty.
\]
Here the constants may depend on the stopping bounds and $\nu$, and both
projection and polarization are bounded at these Sobolev orders.

The first two equations of \eqref{eq:lp-amplitudes} form a wave block
closed under the projected linear evolution, with forcing in both
equations. On the stopped interval, its nondiffusive matrix is bounded
and commutes with the scalar diffusion. Hence
\[
 \|\mathcal S_W(t,r)F\|_{H^{s+1/2}}
 \le C\bigl[1+(\nu(t-r))^{-3/4}\bigr]\|F\|_{H^{s-1}},
  \quad 0\le r<t\le\tau.
\]
The multiplier $\mathcal S(t,r)$ satisfies the same bound. The assumption
$E_{\mathrm{nz}}(0)<\infty$ places the initial wave in $H^{s+1/2}$.
By the Duhamel formula and the integrability of $(t-r)^{-3/4}$, the wave is
continuous in this space up to $t=0$.

We next pair the wave equation at order $s+1/2$ and estimate the forcing
by duality with $H^{s+3/2}$. Young's inequality yields
\[
 \frac{\dd}{\dd t}\|(v,\widehat{\theta_{\mathrm{nz}}^\sh})\|_{H^{s+1/2}}^2
 +c\nu\|(v,\widehat{\theta_{\mathrm{nz}}^\sh})\|_{H^{s+3/2}}^2
 \le C\|(v,\widehat{\theta_{\mathrm{nz}}^\sh})\|_{H^{s+1/2}}^2
       +C\nu^{-1}\|\mathcal F\|_{H^{s-1/2}}^2.
\]
Here $C$ includes the inhomogeneous low-frequency term. Integration gives
the asserted wave regularity. The regularity of the passive
and zero modes follows from the local $H^s$ bounds and the boundedness
of the polarization and pullback maps.

To justify the energy identities, we introduce reference Fourier cutoffs.
The wave nonlinear pairings
converge in $L_t^1$ by $L_t^2H^{s-1/2}$--$L_t^2H^{s+3/2}$ duality,
and the passive and zero-mode pairings converge by
$L_t^2H^{s-1}$--$L_t^2H^{s+1}$ duality. The same regularity gives
convergence of the diffusion terms. The group amplitudes have finite
variation, so It\^o's formula acts only on the frequency weights
and symmetrizers, whose normalized first two Brownian derivatives are
bounded on the stopped interval. After stopping the polarization
energy, dominated convergence and the It\^o isometry give convergence
of the martingale terms. We thereby obtain the combined nonlinear It\^o identities without
assuming that the noise operator maps $H^s$ to $H^{s-1}$.
Letting the stopping bounds increase proves the assertions for every finite
$\tau<\tau_*$.
\end{proof}

\section{Linear enhanced dissipation}\label{sec:lp-linear}

We now turn to the linear decay estimates used in the nonlinear argument.
The simple-zero modes are mixed by the random spanwise shear, whereas
the fast modes are also transported by the Couette shear. We estimate both evolutions
conditionally on the data at the starting time.
Throughout this section, $\nu\in(0,1]$, $\rho\ge0$, and $\beta>1/2$.
Let $W$ be a real Brownian motion relative to a filtration $(\cF_t)$
satisfying the usual conditions. Fourier transforms are normalized by
Plancherel's identity. Constants $C_\beta$ may diverge as $\beta\downarrow1/2$.

\subsection{Characteristics and polarization}

For $(k,\ell)\ne(0,0)$, set $r(\eta)=(h^2+\eta^2)^{1/2}$
and write $\widehat u=e_gv+e_vb$ in the basis
\eqref{eq:intro-polarization}. Along the Stratonovich characteristic
\begin{equation}\label{eq:lp-char}
 \eta_t=\eta_s-k(t-s)-\rho\ell(W_t-W_s),
  \quad r_t=(h^2+\eta_t^2)^{1/2},
\end{equation}
the linearized Fourier equations become
\begin{equation}\label{eq:lp-amplitudes}
 \dot v=\left(-\nu r_t^2+\frac{k\eta_t}{r_t^2}\right)v
                  -\frac{\beta h}{r_t}\theta, \quad
 \dot\theta=-\nu r_t^2\theta+\frac{\beta h}{r_t}v, \quad
 \dot b=-\nu r_t^2b-\frac{\ell}{r_t}v.
\end{equation}
The amplitudes have finite variation. The passive combination satisfies
\begin{equation}\label{eq:lp-passive}
 p=b+\frac{\ell}{\beta h}\theta, \quad
 p(t)=\exp\left(-\nu\int_s^t r_a^2\dd a\right)p(s).
\end{equation}
The coefficients are continuous, and $|k\eta|/r^2\le|k|/(2h)\le1/2$.
It follows that the characteristic system has a pathwise solution on
every finite interval. The estimates below allow us to extend this
solution by density to the stated spaces.

To derive the amplitude equations, we use $\partial_\eta e_g=-\frac{h(k,\eta,\ell)}{r^3}$,
$\partial_\eta e_v=0$, and $\mathcal P(k,\eta,\ell)\partial_\eta e_g=0$,
so $\mathcal P\partial_\eta(e_gv+e_vb)=e_g\partial_\eta v+e_v\partial_\eta b$.
This identity applies to $\mathcal P(y\partial_z)$ and its iterated square.
The shear contributes $k\eta v/r^2$ and $-\ell v/r$, while buoyancy gives
the off-diagonal terms in \eqref{eq:lp-amplitudes}.

\subsection{A Gaussian kernel with drift}

The conditional damping estimates reduce to a Gaussian kernel for the
integrated squared frequency. We first compute this kernel in the
presence of the deterministic drift.

\begin{lemma}\label{lem:lp-gaussian}
Let $a\ge0$, $q>0$, $k,x\in\R$, and let $B$ be a standard real Brownian motion. Define
\[
 G_{q,k,a}(T,x)
 =\E\exp\left[-q\int_0^T(x-k\tau-aB_\tau)^2\dd\tau\right].
\]
For $a>0$, put $\alpha=a\sqrt{2q}$. Then
\begin{equation*}
 \begin{split}
 G_{q,k,a}(T,x)
 ={}&[\cosh(\alpha T)]^{-1/2}\\
 &\times\exp\left[\begin{aligned}
 &-\frac{\alpha}{2a^2}\tanh(\alpha T)x^2
 +\frac{k}{a^2}\bigl(1-\operatorname{sech}(\alpha T)\bigr)x -\frac{k^2}{2a^2}
       \left(T-\frac{\tanh(\alpha T)}{\alpha}\right)
 \end{aligned}\right].
 \end{split}
\end{equation*}
Taking the supremum in the initial frequency gives
\begin{equation}\label{eq:lp-gaussian-sup}
 \sup_{x\in\R}G_{q,k,a}(T,x)
 =[\cosh(\alpha T)]^{-1/2}
  \exp\left[-\frac{k^2}{2a^2}
      \left(T-\frac{2}{\alpha}\tanh\frac{\alpha T}{2}\right)\right].
\end{equation}
For $a=0$, the formula has the continuous limit
\begin{equation}\label{eq:lp-gaussian-zero}
 G_{q,k,0}(T,x)
  =\exp\left[-qT\left(x-\frac{kT}{2}\right)^2
                  -\frac{qk^2T^3}{12}\right].
\end{equation}
\end{lemma}

\Needspace{4\baselineskip}
\begin{proof}
The expectation satisfies the Feynman--Kac equation
$\partial_TG=\frac{a^2}{2}\partial_{xx}G-k\partial_xG-qx^2G$,
with $G(0,x)=1$. Substituting $G=A(T)\exp[-b(T)x^2+d(T)x]$
and comparing coefficients gives $b'=q-2a^2b^2$, $d'=2kb-2a^2bd$,
and $\frac{A'}A=\frac{a^2d^2}{2}-a^2b-kd$.
Solving the system with $b(0)=d(0)=0$ and $A(0)=1$ gives
\[
 \begin{gathered}
 b=\frac{\alpha}{2a^2}\tanh(\alpha T),\quad
 d=\frac{k}{a^2}(1-\operatorname{sech}(\alpha T)),\\
 \log A=-\frac12\log\cosh(\alpha T)
 -\frac{k^2}{2a^2}\left(T-\frac{\tanh(\alpha T)}\alpha\right).
 \end{gathered}
\]
To verify the representation of the expectation, we apply It\^o's formula
to $G(T-\tau,x-k\tau-aB_\tau)$ multiplied by the killing factor,
and stop at $|x-k\tau-aB_\tau|=n$.
Completing the square with
$\tanh z+(1-\operatorname{sech}z)^2/\tanh z=2\tanh(z/2)$
gives \eqref{eq:lp-gaussian-sup} and bounds the candidate by one.
We can therefore remove the stopping by dominated convergence and
recover the expectation. Direct integration gives \eqref{eq:lp-gaussian-zero}.
\end{proof}

\begin{proposition}\label{prop:lp-passive-joint}
Let $k\ne0$, $\ell\in\mathbb Z$, $a=\rho|\ell|$, and
$\lambda_k=\nu^{1/3}|k|^{2/3}$. For $\cF_s$-measurable characteristic initial data $p_s$, the conditional second moment is given exactly by
\begin{equation}\label{eq:lp-p-kernel}
 \E[|p(t)|^2\mid\cF_s]
 =e^{-2\nu h^2T}G_{2\nu,k,a}(T,\eta_s)|p_s|^2,
  \quad T=t-s.
\end{equation}
For some numerical constants $c,C>0$, independent of all parameters,
\begin{equation*}
 \E[|p(t)|^2\mid\cF_s]
 \le C e^{-c(\nu h^2+\lambda_k+\rho|\ell|\sqrt\nu)T}|p_s|^2.
\end{equation*}
The estimate holds for every $\rho\ge0$ and remains valid after
integration against any nonnegative $\cF_s$-measurable reference-frequency
weight. In particular, it applies to the squared reference-frequency
Sobolev norm of a fixed $(k,\ell)$ mode.
\end{proposition}

\begin{proof}
By independence of future Brownian increments, \eqref{eq:lp-passive}
gives \eqref{eq:lp-p-kernel}. For $a>0$, \eqref{eq:lp-gaussian-sup} yields
\begin{equation*}
 \begin{split}
 \sup_x e^{-2\nu h^2T}G_{2\nu,k,a}(T,x)
 ={}&e^{-2\nu h^2T}[\cosh(2a\sqrt\nu T)]^{-1/2}\exp\left[-\frac{k^2}{2a^2}
        \left(T-\frac{\tanh(a\sqrt\nu T)}{a\sqrt\nu}\right)\right].
 \end{split}
\end{equation*}
Write $\mu_a=a\sqrt\nu$. For some numerical $c_0>0$,
$1-\tanh z/z\ge c_0\min\{z^2,1\}$ for $z\ge0$.
For $z\le1$, this follows by integrating
$\tanh^2y\ge(\tanh1)^2y^2$; for $z\ge1$, it follows by monotonicity.
If $\mu_a\le\lambda_k$ and $T\ge\lambda_k^{-1}$, the last exponent
has magnitude at least
\[
 \frac{c_0\lambda_k^3T}{2\mu_a^2}
                \min\{\mu_a^2T^2,1\}
 \ge \frac{c_0}{2}\lambda_kT.
\]
Enlarging the constant extends this bound to $T<\lambda_k^{-1}$.
The hyperbolic-cosine factor is at most $\sqrt2e^{-\mu_aT}$ and
therefore also handles the case $\mu_a\ge\lambda_k$. For $a=0$, use \eqref{eq:lp-gaussian-zero}
and split at $T=\lambda_k^{-1}$. Conditional Tonelli yields the
weighted estimate.
\end{proof}

\Needspace{4\baselineskip}
\begin{remark*}
The passive estimate for $p$ allows arbitrary noise amplitudes. For
the coupled wave system, the estimate below requires small noise.
\end{remark*}

\subsection{Energy and moments of simple-zero modes}

For the simple-zero modes, the wave coupling is skew-adjoint along each
characteristic. This gives an exact energy identity, from which we
derive conditional moment bounds.

For $k=0$, $\ell\ne0$, let
$Q=\widehat u^1+\widehat\theta/\beta=\operatorname{sgn}(\ell)p$.
Along $\eta_t=\eta_s-\rho\ell(W_t-W_s)$, equations
\eqref{eq:lp-amplitudes} give
\begin{equation}\label{eq:lp-zero-exact}
 (v\pm i\theta)(t)=\exp\left[-D_{s,t}
       \ \pm i\beta|\ell|\int_s^t r_a^{-1}\dd a\right]
       (v\pm i\theta)(s), \quad Q(t)=e^{-D_{s,t}}Q(s),
\end{equation}
where $D_{s,t}=\nu\int_s^t(\ell^2+\eta_a^2)\dd a$. In particular,
\begin{equation}\label{eq:lp-zero-density}
 (|v|^2+|\theta|^2+|Q|^2)(t)
 =e^{-2D_{s,t}}(|v|^2+|\theta|^2+|Q|^2)(s).
\end{equation}
With the global reference label $\zeta=\eta+\rho\ell W_t$, we define,
for $\sigma\ge0$,
\begin{equation*}
 \begin{split}
 \mathcal Z_\sigma(t)
 &=\sum_{\ell\ne0}\int_\R\av{\zeta,\ell}^{2\sigma}
       (|v|^2+|\theta|^2+|Q|^2)
              (t,\zeta-\rho\ell W_t,\ell)\dd\zeta,\\
 \mathcal D^0_\sigma(t)
 &=\sum_{\ell\ne0}\int_\R\av{\zeta,\ell}^{2\sigma}
       [\ell^2+(\zeta-\rho\ell W_t)^2]
       (|v|^2+|\theta|^2+|Q|^2)
              (t,\zeta-\rho\ell W_t,\ell)\dd\zeta .
 \end{split}
\end{equation*}
Orthogonality and $u^1=Q-\theta/\beta$ imply that $\mathcal Z_\sigma$
is equivalent to the squared adapted $H^\sigma$ norm of
$\widetilde{(u,\theta)_0}$, with equivalence constants depending only on $\beta$.

\begin{proposition}
\label{prop:lp-zero-moments}
For every $r\ge1$, $\sigma\ge0$, and $0\le s\le t$,
\begin{equation}\label{eq:lp-zero-moments}
 \E[\mathcal Z_\sigma(t)^r\mid\cF_s]
 \le \sqrt2\,e^{-(2r\nu+\sqrt r\,\rho\sqrt\nu)(t-s)}
                         \mathcal Z_\sigma(s)^r.
\end{equation}
No smallness condition on $\rho$ is needed. Almost surely,
\begin{equation}\label{eq:lp-zero-dissipation}
 \mathcal Z_\sigma(t)+2\nu\int_s^t\mathcal D^0_\sigma(a)\dd a
       =\mathcal Z_\sigma(s).
\end{equation}
For $\rho>0$, the first-moment decay rate is sharp on every fixed
$\ell\ne0$. More precisely, for deterministic nonzero initial data on that
mode, with energy density $f(\eta)\ge0$, $f\in L^1(\R)$,
\begin{equation}\label{eq:lp-zero-sharp}
 \begin{split}
 \lim_{t\to\infty}
 e^{(2\nu\ell^2+\rho|\ell|\sqrt\nu)t}
 \E\int_\R(|v|^2+|\theta|^2+|Q|^2)(t,\eta,\ell)\dd\eta=\sqrt2\int_\R
        e^{-\sqrt\nu\,\eta^2/(\rho|\ell|)}f(\eta)\dd\eta>0.
 \end{split}
\end{equation}
\end{proposition}

\begin{proof}
Integration of \eqref{eq:lp-zero-density} in $\zeta$ gives
\eqref{eq:lp-zero-dissipation}. To estimate the moments, we apply
Lemma \ref{lem:lp-gaussian} with
$k=0$, $q=2r\nu$, and $a=\rho|\ell|$ to obtain
\begin{equation*}
\begin{split}
 \E[e^{-2rD_{s,t}}\mid\cF_s]
 ={}&e^{-2r\nu\ell^2T}
 [\cosh(2\sqrt r\,\rho|\ell|\sqrt\nu T)]^{-1/2}\exp\left[-\frac{\sqrt{r\nu}}{\rho|\ell|}
      \tanh(2\sqrt r\,\rho|\ell|\sqrt\nu T)\eta_s^2\right]
 \end{split}
\end{equation*}
for $\rho>0$, with the continuous deterministic limit for $\rho=0$.
Consequently,
\[
 \left(\E[e^{-2rD_{s,t}}\mid\cF_s]\right)^{1/r}
 \le 2^{1/(2r)}e^{-(2\nu\ell^2+\rho|\ell|\sqrt\nu/\sqrt r)T}.
\]
Conditional Minkowski's inequality and $|\ell|\ge1$ yield
\eqref{eq:lp-zero-moments}. We remove the frequency truncations by
monotone convergence. To prove \eqref{eq:lp-zero-sharp}, take $r=1$
and multiply by the exponential factor in that formula. The integrand converges to
$\sqrt2 e^{-\sqrt\nu\,\eta^2/(\rho|\ell|)}f(\eta)$ and is bounded by $\sqrt2f$.
Dominated convergence yields the finite positive limit, proving
sharpness for these data.
\end{proof}

\subsection{Conditional estimates for simple-zero modes}

The nonlinear argument also requires estimates for adapted sources.
We first bound the homogeneous simple-zero evolution in Wiener norms
and estimate its physical gradient after integration in time. Conditioning at
each forcing time then gives the corresponding bounds for the
inhomogeneous equation.
For a simple-zero field $f=(f_u,f_\theta)$ with $\nabla\cdot f_u=0$,
define the adapted and physical Wiener norms, respectively, by
the following formulas, using the Euclidean norm for Fourier vectors:
\begin{equation*}
\begin{split}
 \norm f_{\cA^\sigma_\ad(t)}
 &=\sum_{\ell\ne0}\int_\R
       \av{\eta+\rho\ell W_t,\ell}^{\sigma}
                          |\widehat f(0,\eta,\ell)|\dd\eta,\\
 \norm f_{\cA^1}
 &=\sum_{\ell\ne0}\int_\R
       \av{\eta,\ell}|\widehat f(0,\eta,\ell)|\dd\eta .
 \end{split}
\end{equation*}
Let $S_0(t,s)$ be the evolution \eqref{eq:lp-zero-exact}.
For $\mu=\nu+\rho\sqrt\nu$, set
$\mathcal M=\mu^{-1}+\rho\mu^{-3/2}$.

\begin{proposition}
\label{prop:lp-zero-Wiener}
Let $f_s$ be $\cF_s$-measurable. For every $\sigma\ge0$,
\begin{align}
 \E[\norm{S_0(t,s)f_s}_{\cA^\sigma_\ad(t)}\mid\cF_s]
 &\le C_\beta e^{-c\mu(t-s)}
                    \norm{f_s}_{\cA^\sigma_\ad(s)},\label{eq:lp-Wiener}\\
 \E[\norm{\nabla S_0(t,s)f_s}_\infty\mid\cF_s]
 &\le C_\beta(1+\rho\sqrt{t-s})e^{-c\mu(t-s)}
                    \norm{f_s}_{\cA^1},
                    \label{eq:lp-gradient}\\
 \E\left[\int_s^\infty
          \norm{\nabla S_0(t,s)f_s}_\infty\dd t
                  \,\middle|\,\cF_s\right]
 &\le C_\beta\mathcal M
                    \norm{f_s}_{\cA^1}.
                    \label{eq:lp-gradient-integral}
\end{align}
Here $c>0$ is numerical, and the bounds are uniform in $s$ and
physical frequency.
\end{proposition}

\Needspace{6\baselineskip}
\begin{proof}
The map between $(u,\theta)$ and $(v,\theta,Q)$ is bounded by
$C_\beta$, uniformly in frequency. Thus
$|\widehat{S_0(t,s)f_s}(\eta_t,\ell)|
\le C_\beta e^{-D_{s,t}}|\widehat f_s(\eta_s,\ell)|$.
The reference label $\eta_t+\rho\ell W_t=\eta_s+\rho\ell W_s$ remains
constant along the characteristic. Lemma \ref{lem:lp-gaussian} gives
$\E[e^{-D_{s,t}}\mid\cF_s]
\le\sqrt2e^{-(\nu\ell^2+\rho|\ell|\sqrt{\nu/2})T}$.
Conditional Tonelli's theorem gives \eqref{eq:lp-Wiener}.
To estimate the physical gradient, we use $r_t\le r_s+\rho|\ell||W_t-W_s|$.
Combining Cauchy--Schwarz with the exact second-moment damping, we obtain
\begin{equation*}
\begin{split}
 \E[r_te^{-D_{s,t}}\mid\cF_s]
 &\le\sqrt2r_s e^{-(\nu\ell^2+\rho|\ell|\sqrt{\nu/2})T}+2^{1/4}\rho|\ell|\sqrt T\,
           e^{-(\nu\ell^2+\rho|\ell|\sqrt\nu/2)T}\\
 &\le C r_s(1+\rho\sqrt T)e^{-c\mu T}.
 \end{split}
\end{equation*}
Fourier inversion and $r_s\le\av{\eta_s,\ell}$ yield
\eqref{eq:lp-gradient}. The time-integrated estimate \eqref{eq:lp-gradient-integral} follows by
integration in $T$, since
$\int_0^\infty(1+\rho\sqrt T)e^{-c\mu T}\dd T
\le C(\mu^{-1}+\rho\mu^{-3/2})$.
\end{proof}

\begin{proposition}\label{prop:lp-zero-forced}
Let $F_a$ be progressively measurable, supported on simple-zero modes,
with divergence-free velocity part. Then the mild solution
$U(t)=S_0(t,s)f_s+\int_s^tS_0(t,a)F_a\dd a$
satisfies the estimates below whenever their right-hand sides are finite:
\begin{align}
 \E[\norm{U(t)}_{\cA^\sigma_\ad(t)}\mid\cF_s]
 &\le C_\beta e^{-c\mu(t-s)}
                 \norm{f_s}_{\cA^\sigma_\ad(s)}\notag\\*
 &\quad+C_\beta\int_s^t e^{-c\mu(t-a)}
       \E[\norm{F_a}_{\cA^\sigma_\ad(a)}\mid\cF_s]\dd a,
       \label{eq:lp-forced-Wiener}\\
 \E\left[\int_s^\infty\norm{U(t)}_{\cA^\sigma_\ad(t)}
                       \dd t\,\middle|\,\cF_s\right]
 &\le C_\beta\mu^{-1}
  \Biggl(\norm{f_s}_{\cA^\sigma_\ad(s)} +\E\left[\int_s^\infty\norm{F_a}_{\cA^\sigma_\ad(a)}
                       \dd a\,\middle|\,\cF_s\right]\Biggr),
       \label{eq:lp-forced-Wiener-integral}\\
 \E\left[\int_s^\infty\norm{\nabla U(t)}_\infty
                       \dd t\,\middle|\,\cF_s\right]
 &\le C_\beta\mathcal M
  \Biggl(\norm{f_s}_{\cA^1}\ +\E\left[\int_s^\infty\norm{F_a}_{\cA^1}
                       \dd a\,\middle|\,\cF_s\right]\Biggr).
       \label{eq:lp-forced-gradient}
\end{align}
\end{proposition}

\begin{proof}
The variation-of-constants formula is a pathwise Lebesgue integral.
At each forcing time, we apply Proposition \ref{prop:lp-zero-Wiener}
conditionally on $\cF_a$. By the tower property, the triangle inequality,
and conditional Tonelli's theorem, we obtain
\eqref{eq:lp-forced-Wiener}. Integration in $t$ yields
\eqref{eq:lp-forced-Wiener-integral}. Applying
\eqref{eq:lp-gradient-integral} to $F_a$ and interchanging the
nonnegative integrals gives \eqref{eq:lp-forced-gradient}.
Finally, time and frequency truncations extend the bounds from bounded
smooth sources to general sources.
\end{proof}

\begin{remark*}
The source norm in \eqref{eq:lp-forced-gradient} is physical, whereas
\eqref{eq:lp-forced-Wiener} uses the adapted weight. Because the shift
$\rho\ell W_a$ prevents a uniform comparison in $a$, nonlinear fluxes
require a different estimate. We keep the physical divergence outside
the flux and use Lemma~\ref{lem-smooth}.
For deterministic simple-zero fields, Cauchy--Schwarz gives
$\norm f_{\cA^1}\le C_q\norm f_{H^q}$ for $q>2$.
\end{remark*}

\subsection{A symmetrizer for fast waves}

For the fast modes, the shear term remains in the wave block.
For $k\ne0$, introduce the rescaled amplitudes and symmetrizer
\[
 \begin{aligned}
 X&=\sqrt r\binom v\theta,\quad
 A=\begin{pmatrix}k\eta/(2r^2)&-\beta h/r\\
                   \beta h/r&-k\eta/(2r^2)\end{pmatrix},\\
 c(\eta)&=\frac{k\eta}{2\beta hr},\quad
 M(\eta)=\begin{pmatrix}1&-c(\eta)\\-c(\eta)&1\end{pmatrix}.
 \end{aligned}
\]
Then $A^*M+MA=0$ and
$(1-(2\beta)^{-1})I\le M\le(1+(2\beta)^{-1})I$,
with $1-(2\beta)^{-1}>0$.
Choose $K>(2\beta-1)^{-1}$ and set
\begin{equation*}
 \begin{aligned}
 w(\eta)=\exp\left[K\frac{k\eta}{hr}\right], \quad
 A_k=(|k|/\nu)^{1/3},\quad
 m(\eta)=\begin{cases}
   \exp[\operatorname{sgn}(k)\arctan(\eta/A_k)],&h<A_k,\\
   1,&h\ge A_k.
 \end{cases}
 \end{aligned}
\end{equation*}
Let $\mathsf N=wmM$, $F=X^*\mathsf N X$, and $F_p=m|p|^2$.
All scalar weights and their inverses are uniformly bounded.

\begin{lemma}\label{lem:lp-fast-Ito}
Along \eqref{eq:lp-char}, the wave energy satisfies
\begin{equation}\label{eq:lp-F-Ito}
 \begin{split}
 \dd F={}&\left[-2\nu r^2F-kX^*\mathsf N'X
   +\frac{\rho^2\ell^2}{2}X^*
        \left(\mathsf N''+\frac{2\eta}{r^2}\mathsf N'
                      +\frac{h^2}{r^4}\mathsf N\right)X\right]\dd t\\
 &\quad+\rho X^*
          \left(-\ell \mathsf N'-\frac{\ell\eta}{r^2}\mathsf N\right)X\dd W_t .
 \end{split}
\end{equation}
The passive energy satisfies
\begin{equation}\label{eq:lp-Fp-Ito}
 \dd F_p=
  \left[-2\nu r^2-\frac{km'}m+
                    \frac{\rho^2\ell^2m''}{2m}\right]F_p\dd t
             -\rho\ell\frac{m'}mF_p\dd W_t .
\end{equation}
In particular, each $H=F$ or $H=F_p$ satisfies
\begin{equation}\label{eq:lp-fast-drift}
 \dd H\le[-\nu r_t^2-\lambda_k+C_\beta\rho^2]H\dd t
                      +\rho g_tH\dd W_t, \quad
 |g_t|\le C_\beta,\quad \lambda_k=\nu^{1/3}|k|^{2/3}.
\end{equation}
At zeros of $H$, take $g_t=0$.
\end{lemma}

\begin{proof}
The rescaled amplitude satisfies the Stratonovich equation
\[
 \dd X=(-\nu r^2I+A)X\dd t-\frac{\rho\ell\eta}{2r^2}X\circ\dd W_t.
\]
The It\^o correction is $\rho^2\ell^2(2h^2-\eta^2)X/(8r^4)$.
For $(\eta,X)$, the noise vector field divided by $\rho$ is
\[
 {\cal V}=-\ell\partial_\eta-\frac{\ell\eta}{2r^2}
            \left(X\cdot\partial_X+
                     \overline X\cdot\partial_{\overline X}\right).
\]
Complex derivatives here act on the real and imaginary parts of $X$. Hence
$\mathcal V F=X^*(-\ell\mathsf N'-\ell\eta\mathsf N/r^2)X$ and
\[
 {\cal V}^2F
 =\ell^2X^*
       \left(\mathsf N''+\frac{2\eta}{r^2}\mathsf N'+\frac{h^2}{r^4}\mathsf N\right)X .
\]
Since $A^*\mathsf N+\mathsf N A=0$, this proves
\eqref{eq:lp-F-Ito}. Applying the scalar It\^o formula to
\eqref{eq:lp-passive} gives \eqref{eq:lp-Fp-Ito}.

Since $c'=kh/(2\beta r^3)$,
$|kX^*M'X|\le\frac{k^2h}{(2\beta-1)r^3}X^*MX$.
By the choice of the weight, $-kw'X^*MX=-Kk^2h\,wX^*MX/r^3$
absorbs $-kwX^*M'X$. For $h<A_k$,
$km'/m=\lambda_k/[1+(\eta/A_k)^2]$ and
\[
 \nu r^2+\frac{km'}m
 \ge\lambda_k\left[
       (\eta/A_k)^2+\frac1{1+(\eta/A_k)^2}\right]\ge\lambda_k .
\]
For $h\ge A_k$, we use $\nu r^2\ge\nu h^2\ge\lambda_k$.
In both cases, the deterministic drift is bounded above by
$-(\nu r^2+\lambda_k)H$.
It remains to control the noise terms. Differentiation gives
$\|\ell \mathsf N'\|+\|\ell^2\mathsf N''\|\le C_\beta$ and
$\left|\ell\frac{m'}m\right|+
\left|\ell^2\frac{m''}m\right|\le C_\beta$.
Indeed, the first two derivatives of $M,w$ are bounded, respectively,
by multiples of $h^{-1},h^{-2}$, and $|\ell|\le h$.
The corresponding derivatives of $m$ are bounded by multiples
of $A_k^{-1},A_k^{-2}$ in the region where $|\ell|<A_k$.
Using $|\ell\eta|/r^2\le1/2$, $\ell^2h^2/r^4\le1$, and coercivity
of $\mathsf N$, we bound the It\^o correction by $C_\beta\rho^2H$
and the stochastic coefficient by $C_\beta\rho H$.
This proves \eqref{eq:lp-fast-drift}.
\end{proof}

\begin{proposition}
\label{prop:lp-fast-modes}
For every $r\ge1$ there is $\varepsilon_{\beta,r}>0$ such that,
if $\rho^2\le\varepsilon_{\beta,r}\lambda_k$, then each $H=F$
or $H=F_p$ satisfies
\begin{equation*}
 \E[H(t)^r\mid\cF_s]\le
       e^{-r\lambda_k(t-s)/2}H(s)^r .
\end{equation*}
In particular, $\rho^2\le\varepsilon_{\beta,r}\nu^{1/3}$ suffices for all $k\ne0$.
\end{proposition}

\begin{proof}
It\^o's formula and \eqref{eq:lp-fast-drift} give
\[
 \dd H^r\le[-r\nu r_t^2-r\lambda_k+
                C_{\beta,r}\rho^2]H^r\dd t
                    +r\rho g_tH^r\dd W_t .
\]
For nonintegral $r$, we smooth $x^r$ at zero and pass to the limit.
We choose $\varepsilon_{\beta,r}$ so that
$C_{\beta,r}\rho^2\le r\lambda_k/2$. Then
$e^{r\lambda_k(t-s)/2}H(t)^r$ is a nonnegative local supermartingale.
The asserted bound follows by localization and conditional Fatou's lemma.
\end{proof}

\subsection{Conditional maximal and dissipation estimates}

We next integrate the modewise energies over the reference frequencies
to obtain moment, maximal, and dissipation estimates for the fast field.
Set
$\zeta=\eta+kt+\rho\ell W_t$ and
$\eta_t=\zeta-kt-\rho\ell W_t$. With amplitudes and multipliers
evaluated at $(t,k,\eta_t,\ell)$, define
\begin{equation}\label{eq:lp-fast-field-energy}
 \begin{split}
 \mathcal E_\sigma(t)
 &=\sum_{k\ne0,\ell}\int_\R
         \av{k,\zeta,\ell}^{2\sigma}(F+F_p)\dd\zeta,\quad
 \mathcal D_\sigma(t)
 =\sum_{k\ne0,\ell}\int_\R
         \av{k,\zeta,\ell}^{2\sigma}r_t^2(F+F_p)\dd\zeta,\\
 \mathcal L_\sigma(t)
 &=\sum_{k\ne0,\ell}\int_\R
         \av{k,\zeta,\ell}^{2\sigma}\lambda_k(F+F_p)\dd\zeta .
 \end{split}
\end{equation}
In particular, $\mathcal L_\sigma\ge\nu^{1/3}\mathcal E_\sigma$.
The energy is equivalent to
\[
 \sum_{k\ne0,\ell}\int_\R\av{k,\zeta,\ell}^{2\sigma}
       [r_t(|v|^2+|\theta|^2)+|p|^2]\dd\zeta .
\]
It controls the squared adapted $H^\sigma$ norm and, at $t=0$, is bounded
by $C_\beta\norm{(u^{\rm in},\theta^{\rm in})}_{H^{\sigma+1/2}}^2$.
At a general starting time $s$, the factor $r_s$ must be retained:
the energy is no longer uniformly equivalent to the reference-frequency
$H^{\sigma+1/2}$ norm.

\begin{proposition}
\label{prop:lp-fast-field}
Let $\sigma\ge0$, $r\ge1$, and let the initial fast field at time $s$ be $\cF_s$-measurable with $\mathcal E_\sigma(s)<\infty$
almost surely. There is $\varepsilon_{\beta,r}>0$ such that,
if $\rho^2\le\varepsilon_{\beta,r}\nu^{1/3}$, the unique characteristic solution satisfies
\begin{align}
 \E[\mathcal E_\sigma(t)^r\mid\cF_s]
 &\le e^{-r\nu^{1/3}(t-s)/2}
                         \mathcal E_\sigma(s)^r,\label{eq:lp-full-moment}\\
 \E\left[\sup_{t\ge s}
       \left(e^{\nu^{1/3}(t-s)/4}\mathcal E_\sigma(t)\right)^r
                       \,\middle|\,\cF_s\right]
 &\le2\mathcal E_\sigma(s)^r,\label{eq:lp-full-maximal}\\
 \E\left[\int_s^\infty e^{\nu^{1/3}(t-s)/4}
         \bigl(\nu\mathcal D_\sigma(t)+\mathcal L_\sigma(t)\bigr)
                        \dd t\,\middle|\,\cF_s\right]
 &\le2\mathcal E_\sigma(s).\label{eq:lp-full-dissipation}
\end{align}
For \eqref{eq:lp-full-maximal}, the smallness constant is that associated
with moment $2r$. All constants are uniform in $s$, the reference label, the number
of modes, and the final time.
\end{proposition}

\begin{proof}
We first work with finitely many $(k,\ell)$ modes, compact
reference-frequency support, and bounded initial energy. Stochastic Fubini and
Lemma \ref{lem:lp-fast-Ito} give
\begin{equation}\label{eq:lp-integrated-drift}
 \begin{aligned}
 \dd\mathcal E_\sigma
 &\le[-\nu\mathcal D_\sigma-\mathcal L_\sigma
                  +C_\beta\rho^2\mathcal E_\sigma]\dd t
                          +\rho G_t\dd W_t,\\
 |G_t|&\le C_\beta\mathcal E_\sigma .
 \end{aligned}
\end{equation}
All modes are driven by the same $W$, so we first bound their stochastic
coefficients and then integrate over frequency. For $q\ge1$, It\^o's formula gives
\begin{equation*}
\begin{split}
 \dd\mathcal E_\sigma^q\le 
   -q\mathcal E_\sigma^{q-1}
        (\nu\mathcal D_\sigma+\mathcal L_\sigma)\dd t
   +C_{\beta,q}\rho^2\mathcal E_\sigma^q\dd t +q\rho\mathcal E_\sigma^{q-1}G_t\dd W_t .
 \end{split}
\end{equation*}
If $C_{\beta,q}\rho^2\le q\nu^{1/3}/2$, then
$e^{q\nu^{1/3}(t-s)/2}\mathcal E_\sigma(t)^q$ is a
nonnegative supermartingale. This proves \eqref{eq:lp-full-moment}.

To obtain the maximal estimate, take $q=2r$ and
$Y_t=e^{\nu^{1/3}(t-s)/4}\mathcal E_\sigma(t)$. Since $Y_t^q$
is a nonnegative supermartingale, conditional optional stopping at its
first crossing of $(R\mathcal E_\sigma(s))^q$ gives, for $R\ge1$,
\[
 \Prob\left(\sup_{t\ge s}Y_t>R\mathcal E_\sigma(s)
 \,\middle|\,\cF_s\right)\le R^{-q}.
\]
We let the stopping horizon tend to infinity and integrate the tail to obtain
\[
 \E[(\sup_{t\ge s}Y_t)^r\mid\cF_s]
 \le\left(1+\frac{r}{q-r}\right)\mathcal E_\sigma(s)^r
 =2\mathcal E_\sigma(s)^r .
\]
On $\{\mathcal E_\sigma(s)=0\}$ the solution vanishes.

For the dissipation estimate, take $C_\beta\rho^2\le\nu^{1/3}/4$ and multiply
\eqref{eq:lp-integrated-drift} by $e^{\nu^{1/3}(t-s)/4}$.
The drift is at most
$-e^{\nu^{1/3}(t-s)/4}
\left(\nu\mathcal D_\sigma+\frac12\mathcal L_\sigma\right)$.
Taking conditional expectations after localization, and then applying
conditional Fatou's lemma, yields \eqref{eq:lp-full-dissipation}.

To remove the truncations, apply the maximal estimate to differences
of truncated solutions. This gives convergence on finite intervals in
the energy space, both in probability and in locally finite conditional
moments. For unbounded initial energy,
we first restrict to $\{\mathcal E_\sigma(s)\le n\}\in\cF_s$.
The characteristic equations identify the limit and ensure uniqueness.
All three bounds pass to the limit by conditional Fatou's lemma, with
weak lower semicontinuity also giving the dissipation estimate.
\end{proof}

\begin{remark*}

If $\rho=\nu^{1/4}$, then $\rho^2/\nu^{1/3}=\nu^{1/6}\to0$.
It follows from Proposition \ref{prop:lp-fast-field} that every
fixed-order conditional moment is controlled for small $\nu$.
\end{remark*}

\Needspace{6\baselineskip}
\section{Nonlinear streamwise-independent dynamics}\label{sec:nonlinear}
\subsection{Global regularity for arbitrary data}

We first study the dynamics on the $x$-independent subspace, which is invariant under \eqref{eq:model}.
Throughout this section, $\beta>0$ and $\rho\ge0$ are arbitrary.
On $\mathcal D=\R_y\times\T_z$, set
$v=(v^y,v^z)=(u^2,u^3)$, $q=u^1+\beta^{-1}\theta$, and $Y=(v,\theta,q)$.
Let $\nabla=(\partial_y,\partial_z)$, write
$\mathcal P_2$ for the two-dimensional Leray projection, and set $B=y\partial_z$.
The subsystem reads
\begin{equation}\label{eq:2d-subsystem}
\begin{cases}
\mathrm dv+\mathcal P_2\big((v\cdot\nabla)v+\beta\theta e_y\big)\,\mathrm dt
       +\rho\mathcal P_2 Bv\circ\mathrm dW_t=\nu\Delta v\,\mathrm dt,\\
\mathrm d\theta+\big(v\cdot\nabla\theta-\beta v^y\big)\,\mathrm dt
       +\rho B\theta\circ\mathrm dW_t=\nu\Delta\theta\,\mathrm dt,\\
\mathrm dq+v\cdot\nabla q\,\mathrm dt
       +\rho Bq\circ\mathrm dW_t=\nu\Delta q\,\mathrm dt,\\
\nabla\cdot v=0.
\end{cases}
\end{equation}
All three equations are driven by the same Brownian motion. Solutions are
understood through the Leray transport group or the weak Stratonovich
formulation against smooth, compactly supported test functions;
velocity test functions are solenoidal. This formulation does not require the
mapping property $y\partial_z:H^j\to H^{j-1}$.

\Needspace{6\baselineskip}
\begin{proposition}\label{prop:2d-global}
Let $\nu>0$, $\rho\ge0$, $\beta>0$, and let $s\ge3$ be an integer.
For every deterministic initial state
$Y_{\rm in}=(v_{\rm in},\theta_{\rm in},q_{\rm in})\in H^s(\mathcal D)$
with $\nabla\cdot v_{\rm in}=0$, the $x$-independent restriction of
\eqref{eq:model} has a unique adapted, pathwise strong solution on
$[0,\infty)$.
For every finite $T$, almost surely,
$Y\in C([0,T];H^s(\mathcal D))\cap L^2(0,T;H^{s+1}(\mathcal D))$.
No smallness condition, mean-zero condition, or relation between $\nu$ and $\rho$ is required.
Moreover,
\begin{equation}\label{eq:2d-E0-exact}
\|Y(t)\|_2^2+2\nu\int_0^t\|\nabla Y(a)\|_2^2\,\mathrm da
=\|Y_{\rm in}\|_2^2, \quad t\ge0,
\end{equation}
almost surely. Squared norms of $Y$ denote sums over its components.
An $H^{s-1}$-valued interpretation of the stochastic integral requires an
additional generator-domain assumption.
\end{proposition}

\begin{proof}[Proof of Proposition~\ref{prop:2d-global}]
We first estimate the noise contribution. For
$\xi=(\eta,\ell)\ne0$, set $r=(\eta^2+\ell^2)^{1/2}$,
$e_\perp=(\ell,-\eta)/r$, and $\widehat v=e_\perp a$.
The zero frequency is negligible in $L^2$, and the noise vanishes for
$\ell=0$. Since $\mathcal P_2(\xi)\partial_\eta e_\perp=0$, the noise acts on each of
$a,\widehat\theta,\widehat q$ through
$+\rho\ell\partial_\eta\circ\mathrm dW_t$.

Set $w_j(\eta,\ell)=(1+\eta^2+\ell^2)^j$ and
\[
E_j=\|Y\|_{H^j}^2
=\sum_{\ell\in\mathbb Z}\int_{\mathbb R}
w_j\bigl(|a|^2+|\widehat\theta|^2+|\widehat q|^2\bigr)\,\mathrm d\eta.
\]
The noise contribution to $\mathrm dE_j$ is
\begin{align*}
&\frac{\rho^2}{2}\sum_{\ell}\int
\ell^2\partial_\eta^2w_j
\bigl(|a|^2+|\widehat\theta|^2+|\widehat q|^2\bigr)\,\mathrm d\eta\,\mathrm dt
 -\rho\sum_{\ell}\int
\ell\partial_\eta w_j
\bigl(|a|^2+|\widehat\theta|^2+|\widehat q|^2\bigr)\,\mathrm d\eta\,\mathrm dW_t.
\end{align*}
The It\^o correction in the amplitude equations is $\rho^2\ell^2\partial_\eta^2/2$.
Integration by parts puts both derivatives on $w_j$.
For every integer $j\ge1$,
$|\ell\partial_\eta w_j|+|\ell^2\partial_\eta^2w_j|\le C_jw_j$.
Thus the noise drift has absolute value at most $C_j\rho^2E_j$,
and its martingale coefficient satisfies $|b_j|\le C_j\rho E_j$.
Both terms vanish for $j=0$. A stopped group approximation extends these
computations from the smooth generator core to local $H^s$ solutions
without $y$-weighted norms.

For local existence, we use the group generated by
$\operatorname{diag}(\mathcal P_2 B,B,B)$. It translates amplitudes in $\eta$,
parallel transports the velocity basis, and is unitary on $L^2$ and
bounded on $H^j$ for each finite parameter. Conjugating by this group gives
a pathwise parabolic equation with symbol
$-\nu[(\eta-\rho\ell W_t)^2+\ell^2]$. This equation is uniformly elliptic whenever
$|W_t|$ is bounded. In reference coordinates, the inverse group consists
of the scalar pullback followed by an orthogonal Fourier change of frame.
The latter has norm one on $H^j$, so the conjugated nonlinear terms involve
these bounded operators and the derivatives
$\nabla_t=(\partial_Y-\rho W_t\partial_Z,\partial_Z)$.
The Sobolev algebra property gives a locally Lipschitz map
$\mathcal N_t:H^s\to H^{s-1}$ with
$\|\mathcal N_t(Y)\|_{H^{s-1}}
 \le C_s(1+\rho|W_t|)\|Y\|_{H^s}^2$.
On stopped finite intervals, the physical and group norms are equivalent,
and the heat propagator gains one derivative with constants depending on
the Brownian bound. These estimates allow us to apply Picard iteration,
stopped at Brownian and Sobolev exit times, to obtain a unique adapted
maximal solution with a physical $H^s$
continuation criterion.

We next derive bounds that exclude finite-time breakdown. Incompressibility
cancels the $L^2$ transport terms, buoyancy cancels between velocity and
temperature, and the $q$ equation has no source. Combining these
cancellations with the noise identity gives \eqref{eq:2d-E0-exact} up to
the maximal lifetime. In particular, with
$E_0^{\rm in}=\|Y_{\rm in}\|_2^2$,
$\int_0^t\|\nabla v(a)\|_2^2\,\mathrm da\le E_0^{\rm in}/(2\nu)$.

For the $H^1$ estimate, the velocity convection term cancels:
with
$\omega=\partial_yv^z-\partial_zv^y$, incompressibility gives
$\|\nabla v\|_2^2=\|\omega\|_2^2$, and
$\int\omega(v\cdot\nabla\omega)=0$.
For $f=\theta$ or $q$, the differentiated transport satisfies
\[
\left|\sum_i\int\partial_if\,
\bigl((\partial_iv)\cdot\nabla f\bigr)\right|
\le \|\nabla v\|_2\|\nabla f\|_4^2.
\]
On $\mathbb R\times\mathbb T$, the inhomogeneous two-dimensional
Gagliardo--Nirenberg inequality gives
$\|\nabla f\|_4^2
\le C\|\nabla f\|_2
\bigl(\|\nabla^2f\|_2+\|\nabla f\|_2\bigr)$.
Applying Young's inequality, we bound the transport contribution by
\[
\frac\nu4\|\nabla^2f\|_2^2
+C\bigl(\nu^{-1}\|\nabla v\|_2^2+\|\nabla v\|_2\bigr)
\|\nabla f\|_2^2.
\]
Combining the transport estimate with the differentiated buoyancy cancellation, and setting
$D_j^2=\|\nabla Y\|_{H^j}^2$, we obtain
\begin{equation}\label{eq:2d-H1-inequality}
\mathrm dE_1+\nu D_1^2\,\mathrm dt
\le A_1(t)E_1\,\mathrm dt+b_1(t)\,\mathrm dW_t,
 \quad |b_1(t)|\le C\rho E_1(t),
\end{equation}
where $A_1(t)=C\bigl[\rho^2+\nu^{-1}\|\nabla v(t)\|_2^2
+\|\nabla v(t)\|_2\bigr]$.
Let $\tau_*$ be the maximal lifetime. By \eqref{eq:2d-E0-exact}, for every finite $T$,
\[
\int_0^{T\wedge\tau_*}A_1(t)\,\mathrm dt
\le C\bigl[\rho^2T+\nu^{-2}E_0^{\rm in}
+T^{1/2}\nu^{-1/2}(E_0^{\rm in})^{1/2}\bigr]<\infty.
\]
Set $\chi_1=b_1/E_1$ when $E_1>0$, and zero when $E_1=0$.
The stochastic integrating factor for
$\mathrm dX=A_1X\,\mathrm dt+\chi_1X\,\mathrm dW_t$ gives
\[
E_1(t)\le E_1(0)
\exp\left(\int_0^t(A_1-\tfrac12\chi_1^2)\,\mathrm da
+\int_0^t\chi_1\,\mathrm dW_a\right).
\]
Since $|\chi_1|\le C\rho$, the right-hand side is continuous at finite
stopped times. Hence
\begin{equation}\label{eq:2d-H1-finite}
\sup_{t<T\wedge\tau_*}E_1(t)<\infty,
 \quad \int_0^{T\wedge\tau_*}D_1(t)^2\,\mathrm dt<\infty
\end{equation}
almost surely, where the second bound follows from
\eqref{eq:2d-H1-inequality} after stopping at finite $H^1$ levels.

For higher derivatives, let $j\ge2$ be an integer and $1\le a\le j$. Sobolev embedding and
interpolation give $\|\nabla^a f\|_4
\le C_j\|f\|_{H^1}^{1-\gamma_a}
\|f\|_{H^{j+1}}^{\gamma_a}$, where $\gamma_a=(a-\tfrac12)/j$.
The highest-order commutator factors have derivative orders
$a$ and $j+1-a$, so the corresponding interpolation exponents sum to one.
Cancelling the advective term by incompressibility and summing over $v,\theta,q$, we obtain
$\left|\langle (v\cdot\nabla)Y,Y\rangle_{H^j}\right|
\le C_j\|Y\|_{H^1}\|Y\|_{H^j}\|Y\|_{H^{j+1}}$.
Interpolation also controls the lower-order terms. The Leray
projection drops out of the pairing, and buoyancy cancels at each order.
Since $\|Y\|_{H^{j+1}}^2=E_j+D_j^2$, Young's inequality and
the preceding noise bound give
\begin{equation*}
\mathrm dE_j+\nu D_j^2\,\mathrm dt
\le C_j\bigl[\rho^2+\|Y\|_{H^1}
+\nu^{-1}\|Y\|_{H^1}^2\bigr]E_j\,\mathrm dt
+b_j(t)\,\mathrm dW_t,
 \quad |b_j(t)|\le C_j\rho E_j(t).
\end{equation*}
The first-order bound in \eqref{eq:2d-H1-finite} makes the drift coefficient integrable
on finite intervals up to $\tau_*$. Applying the integrating-factor argument once more gives
\[
\sup_{t<T\wedge\tau_*}E_s(t)
+\nu\int_0^{T\wedge\tau_*}D_s(t)^2\,\mathrm dt<\infty
\quad\hbox{almost surely}.
\]
The continuation criterion now gives $\tau_*\ge T$. Letting
$T\to\infty$ through the integers proves global existence; the remaining properties follow
from the local construction.
\end{proof}

\subsection{Conditional divergence smoothing}
To estimate the response to nonlinear fluxes, we use the pullbacks
$V,\Theta,Q$ and the quantities $\nabla_t,S,M$ defined at the start of the next subsection.
They satisfy $\nabla_t\cdot V=0$ and $\overline{V^y}=0$; normalized Haar measure
gives $S^2+M^2=\|Y^\sh\|_{H^s}^2$. By Fourier orthogonality, this is also the squared
parallel-transport norm. All products below are formed from the componentwise pullbacks.

\begin{lemma}\label{lem-smooth}
Let $\Pi_{\mathrm{nz},z}$ be the projection onto the nonzero $z$ Fourier
modes. At reference frequency
$(\eta,\ell)$, set $\xi_r=(\eta-\rho\ell W_r,\ell)$ and define
$\widehat{\mathcal P}_r(\eta,\ell)
 =I_2-\frac{\xi_r\otimes\xi_r}{|\xi_r|^2}$ for $\xi_r\ne0$.
Set the symbol to $I_2$ at $\xi_r=0$; the projection acts as
the identity on scalar components. Let $\mathcal S_0(t,r)$
be the homogeneous simple-zero propagator for \eqref{eq:2d-subsystem}.
Let $\mathcal H_r$ be an $\cF_r$-measurable tensor/scalar flux in global
reference coordinates. Put $\mu=\nu+\rho\sqrt\nu$.
For every $p\ge2$ there are $c_p,C_p>0$ such that, for $t>r$,
\begin{align}
\Big(\E[\|\mathcal S_0(t,r)\Pi_{\mathrm{nz},z}\mathcal P_r
 \operatorname{div}_r\mathcal H_r\|_{H^s}^p\mid\cF_r]\Big)^{1/p}
 &\le \frac{C_p}{\sqrt{\nu(t-r)}}e^{-c_p\mu(t-r)}
       \|\mathcal H_r\|_{H^s},\label{cond-smooth}\\
\Big(\E[\|\mathcal S_0(t,r)Z_r\|_{H^s}^p\mid\cF_r]\Big)^{1/p}
 &\le C_pe^{-c_p\mu(t-r)}\|Z_r\|_{H^s}.\label{cond-hom}
\end{align}
The input $Z_r$ is $\cF_r$-measurable, has zero $z$ mean,
and has $\operatorname{div}_r$-free velocity. The projection acts only on
velocity fluxes. All norms use the same reference Fourier weight.
The estimates hold conditionally for adapted inputs; they are not bounds on the expected operator norm.
\end{lemma}

\begin{proof}
Fix a nonzero $z$ frequency $l$ and condition on $\cF_r$.
Let $\tau=t-r$ and let $\xi$ be the physical $y$ frequency at $r$.
The future characteristic is
$\xi_h=\xi-\rho l(W_{r+h}-W_r)$.
The projected transport is orthogonal and the gravity coupling is skew-adjoint,
so only diffusion changes the norm along the characteristic. The initial squared norm is therefore multiplied by
$\exp\left[-2\nu\int_0^\tau(\xi_h^2+l^2)\dd h\right]$.
The scalar $q$ has the same attenuation. The Brownian harmonic-oscillator
formula gives
\begin{align*}
 \E\exp\left[-p\nu\int_0^\tau(\xi_h^2+l^2)\dd h\right]
 &=e^{-p\nu l^2\tau}\cosh(\sqrt{2p}\rho|l|\sqrt\nu\tau)^{-1/2}
       e^{-a_{p,l}(\tau)\xi^2},\\
 a_{p,l}(\tau)&=\frac{\sqrt{p\nu/2}}{\rho|l|}
       \tanh(\sqrt{2p}\rho|l|\sqrt\nu\tau).
\end{align*}
For $\rho=0$, the continuous limit is $a_{p,l}(\tau)=p\nu\tau$.
This gives \eqref{cond-hom}. To estimate the divergence source, we note that the divergence
at the input time contributes a factor at most
$(\xi^2+l^2)^{1/2}$. The Leray projection has norm one.
Using $\sup_x x^2e^{-ax^2}\le C/a$ and
$a_{p,l}(\tau)^{-1}\le C_p\left[(\nu\tau)^{-1}
                  +\frac{\rho|l|}{\sqrt\nu}\right]$,
the squared pointwise $L^p(\Omega)$ bound for fixed $l$ is at most
\[
 C_p\left[(\nu\tau)^{-1}+\frac{\rho|l|}{\sqrt\nu}+l^2\right]
 e^{-c_p\rho|l|\sqrt\nu\tau-2\nu l^2\tau}.
\]
Absorbing the last two terms into the exponentials gives
$C_p(\nu\tau)^{-1}e^{-c_p\mu\tau}$ uniformly for $|l|\ge1$.
Since $p/2\ge1$, the conditional Minkowski inequality yields
\eqref{cond-smooth}--\eqref{cond-hom} in the reference Sobolev norm.
Here the scalar weight commutes with the Fourier isometries, and the
shift at $r$ is fixed under conditioning. Lemma~\ref{lem:lp-gaussian}
applies to all initial frequencies; density extends the bounds from the smooth core.
\end{proof}

\begin{remark*}
The kernel satisfies
$\int_0^\infty \frac{e^{-c\mu t}}{\sqrt{\nu t}}\dd t
       =\frac{C}{\sqrt{\nu\mu}}$.
For $\rho=\kappa_0\nu^{1/4}$ this is at most
$C\kappa_0^{-1/2}\nu^{-7/8}$. By controlling the physical derivative at
the forcing time, conditional smoothing avoids the additional factor
$|W_r|$ arising from a separate estimate in the reference norm.
\end{remark*}

\subsection{Small-data nonlinear enhanced dissipation}

Let $(V,\Theta,Q)=Y^\sh$ and
$\nabla_t=(\partial_Y-\rho W_t\partial_Z,\partial_Z)$.
A bar denotes the normalized $z$ average. Set
$S(t)=\|\widetilde{Y^\sh}\|_{H^s(\mathcal D)}$ and
$M(t)=\|\overline{(V^z,\Theta,Q)}\|_{H^s(\R)}$.
Finite energy and incompressibility give $\overline{V^y}=0$.
To separate the free heat evolution, set
$m=\overline{(V^z,\Theta,Q)}$,
$m_{\rm nl}(t)=m(t)-e^{\nu t\partial_{YY}}m_{\rm in}$, and
$M_{\rm nl}(t)=\|m_{\rm nl}(t)\|_{H^s(\R)}$, where
$m_{\rm in}=\overline{(v_{\rm in}^z,\theta_{\rm in},q_{\rm in})}$.

\begin{proposition}\label{thm-zero}
Let $s\ge3$ be an integer, let $\nu,\beta>0$ and $\rho\ge0$, and suppose that the deterministic
initial data for the $x$-independent restriction of \eqref{eq:model}
belong to $H^s(\mathcal D)$, with
$\operatorname{div}v_{\rm in}=0$. No initial mean condition is imposed.
Set $\eps=\|(v_{\rm in},\theta_{\rm in},q_{\rm in})\|_{H^s}$ and
$\mu=\nu+\rho\sqrt\nu$.
For every $\zeta\in(0,1)$ there are $c_{s,\zeta},C_{s,\zeta}>0$ and a numerical
$a>0$ such that, if
\begin{equation}
 \eps\le c_{s,\zeta}\sqrt{\nu\mu}
          =c_{s,\zeta}\sqrt{\nu(\nu+\rho\sqrt\nu)},\label{zero-small}
\end{equation}
the unique global strong solution satisfies, with probability at least
$1-\zeta$,
\begin{align}
 \sup_{t\ge0}e^{a\mu t}S(t)
 +\sqrt\nu\|e^{a\mu t}\nabla_t
       \widetilde{Y^\sh}\|_{L^2_tH^s}
 &\le C_{s,\zeta}\eps,\label{zero-wave}\\
 \sup_{t\ge0}M_{\rm nl}(t)
       +\sqrt\nu\|\partial_Ym_{\rm nl}\|_{L^2_tH^s}
 &\le \frac{C_{s,\zeta}\eps^2}{\sqrt{\nu\mu}}.\label{zero-mean}
\end{align}
On the same event, the total mean satisfies
\begin{equation}\label{zero-total-mean}
 \sup_{t\ge0}M(t)+\sqrt\nu\|\partial_Ym\|_{L^2_tH^s}
       \le C_s\eps.
\end{equation}
For $\chi(Y,Z,t)=t\widetilde\Theta(Y,Z,t)/\beta$, the same event also gives
\begin{align}
 \sup_{t\ge0}\|\chi(t)\|_{H^s(\mathcal D)}
       &\le C_{s,\beta,\zeta}\eps\mu^{-1},
                         \label{zero-chart-reference}\\
 \sup_{t\ge0}\|\nabla_t\chi(t)\|_\infty
       &\le C_{s,\beta,\zeta}\eps\mu^{-1}
                         (1+\rho\mu^{-1/2}).\label{zero-chart-physical}
\end{align}
On this event, the mean field tends to zero in $H^s(\R)$, and the
complete perturbation tends to zero in the physical $H^s$ norm.
Up to fixed constants, one may take $c_{s,\zeta}=c_s\sqrt \zeta$ and
$C_{s,\zeta}=C_s\zeta^{-1/2}$. The failure probability is $\zeta$;
the probability estimates use only second and fourth moments.
The constants in the decay and mean bounds need not depend on $\beta$.
Reconstructing $u^1=q-\beta^{-1}\theta$ introduces
$\beta$-dependence.
For $\rho=\kappa_0\nu^{1/4}$, the sufficient condition
$\eps\le c_{s,\zeta}\kappa_0^{1/2}\nu^{7/8}$ implies \eqref{zero-small}.
\end{proposition}

\begin{proof}[Proof of Proposition~\ref{thm-zero}]
The assertion is immediate for $\eps=0$. For $\eps>0$,
Proposition~\ref{prop:2d-global} gives a global solution. The group
approximations used in that proposition and in Proposition~\ref{thm:affine-local-strong}
justify the stopped energy estimates below.

For deterministic $R>\eps$, let $\tau=\tau_R$ be the first time
$(S^2+M^2)^{1/2}=R$, with $\tau_R=\infty$ if no exit occurs.
To specify the nonlinear sources, set
\[
 \mathcal F=(V\otimes V,V\Theta,VQ), \quad
 \mathcal F_{\rm mm}
   =(\overline V\otimes\overline V,
     \overline V\,\overline\Theta,\overline V\,\overline Q).
\]
The velocity block uses row divergence, and the other blocks use scalar
divergence. Since $\overline V=(0,\overline{V^z})$ depends only on $Y$,
$\operatorname{div}_t\mathcal F_{\rm mm}=0$. We therefore define
\[
 \mathcal H=\mathcal F-\mathcal F_{\rm mm}, \quad
 \mathcal F_{\rm s}=\Pi_{\mathrm{nz},z}\mathcal H, \quad
 \mathcal F_{\rm m}=\overline{V^y(V^z,\Theta,Q)}.
\]
Each term in $\mathcal H$ contains a simple-zero factor. Before $\tau$,
\begin{align}
 \|\mathcal F_{\rm s}(t)\|_{H^s}+\|\mathcal H(t)\|_{H^s}
 &\le C_s(M(t)+S(t))S(t)\le C_sR S(t),\label{flux-simple}\\
 \|\mathcal F_{\rm m}(t)\|_{H^s(\R)}&\le C_sS(t)^2.\label{flux-mean}
\end{align}
The second bound follows from $\overline{V^y}=0$: the mean flux is
$\overline{\widetilde V^y\widetilde{(V^z,\Theta,Q)}}$.
These flux bounds are independent of $W_t$. Fourier isometry identifies the norms
with those in the polarization energy. We keep the physical derivative outside
the flux and apply Lemma~\ref{lem-smooth}.

After $\tau$, we continue the simple component by the homogeneous linear evolution.
Denote this continuation by $Z^\tau$ and set
$x_p(t)=(\E\|Z^\tau(t)\|_{H^s}^p)^{1/p}$ for $p=2,4$.
The continuation agrees with the original simple component before $\tau$.
On $\{\tau<\infty\}$, orthogonality of the projected transport and
the skew-adjoint buoyancy coupling give, for $t\ge\tau$,
\[
 \|Z^\tau(t)\|_{H^s}^2
 +2\nu\int_\tau^t\|\nabla_r Z^\tau(r)\|_{H^s}^2\dd r
 =\|Z^\tau(\tau)\|_{H^s}^2.
\]
This identity holds for every $\nu,\beta>0$ and $\rho\ge0$.
Consequently, the reference norm is nonincreasing after exit and is
bounded by $R$ for all times. Its weighted moment suprema on finite
intervals are therefore finite. The stopped flux
$1_{r<\tau}\mathcal F_{\rm s}(r)$ is $\cF_r$-measurable, and
\eqref{flux-simple} bounds one factor by $R$. Lemma~\ref{lem-smooth},
the Minkowski inequality, and the mild formula give
\[
 x_p(t)\le C_p\eps e^{-c_p\mu t}
 +C_{s,p}R\int_0^t\frac{e^{-c_p\mu(t-r)}}{\sqrt{\nu(t-r)}}x_p(r)\dd r.
\]
Choose $b>0$ below both exponential constants. The weighted kernel has
$L^1$ norm $C_{s,p}(\nu\mu)^{-1/2}$. Fix an absorption constant $r_s>0$
small enough for $p=2,4$. If $R\le r_s\sqrt{\nu\mu}$, then, for every finite $T$,
\[
 X_p(T):=\sup_{0\le t\le T}e^{b\mu t}x_p(t)
 \le C_p\eps+\frac12X_p(T).
\]
The preceding continuation bound makes $X_p(T)$ finite. Absorption
therefore gives a bound independent of $T$, and hence
\begin{equation}
 x_p(t)\le C_s\eps e^{-b\mu t}, \quad p=2,4.\label{moment-decay}
\end{equation}
The fourth-moment inequality is linear because one flux factor is bounded by $R$.

We use this moment decay to control the exit probability. The gravity
and stochastic parallel-transport terms cancel in the reference
$H^s$ energy. Pairing the flux with the physical gradient and using
the self-adjoint Leray projection, the commuting scalar multiplier,
and Young's inequality gives the following bound. We first integrate
up to $T\wedge\tau$ for finite $T$. The bounds are independent of $T$,
so monotone convergence then permits $T\to\infty$, both here and in the
weighted energy and mean estimates below:
\[
 \sup_{t\le\tau}(S(t)^2+M(t)^2)
 \le\eps^2+C_s\nu^{-1}R^2\int_0^\tau S(t)^2\dd t.
\]
The reference energy has no martingale term. By \eqref{moment-decay},
\begin{align*}
 \E\sup_{t\le\tau}(S(t)^2+M(t)^2)
 &\le C_s\eps^2\left(1+\frac{R^2}{\nu\mu}\right)
 \le C_s\eps^2,\\
 \mathbb P(\tau_R<\infty)&\le C_s\eps^2/R^2,
\end{align*}
where continuity at the exit time gives the second bound.

For the weighted simple-component estimate, set
$D=\|\nabla_t\widetilde{Y^\sh}\|_{H^s}$. The same argument gives
$\dd S^2+\nu D^2\dd t
       \le C_s\nu^{-1}\|\mathcal F_{\rm s}\|_{H^s}^2\dd t$.
Multiply by $e^{2a\mu t}$ with $0<a<b$.  Using
\eqref{flux-simple} and \eqref{moment-decay} in the integrated inequality
gives
\begin{equation}
 \E\left[\sup_{t\le\tau}e^{2a\mu t}S(t)^2
       +\nu\int_0^\tau e^{2a\mu t}D(t)^2\dd t\right]
 \le C_s\eps^2\left(1+\frac{R^2}{\nu\mu}\right)
 \le C_s\eps^2.\label{simple-max}
\end{equation}
Differentiating the weight contributes
$2a\mu\int e^{2a\mu t}S(t)^2\dd t$, which is bounded by the moment estimate.

For the mean response, we use the equation
$(\partial_t-\nu\partial_{YY})m=-\partial_Y\mathcal F_{\rm m}$.
The response $m_{\rm nl}$ starts from zero even when
$m_{\rm in}\ne0$. Combining its energy estimate with
\eqref{flux-mean} and \eqref{moment-decay} for $p=4$, we obtain
\begin{equation}
 \E\left[\sup_{t\le\tau}M_{\rm nl}(t)^2
       +\nu\int_0^\tau\|\partial_Ym_{\rm nl}\|_{H^s}^2\dd t\right]
 \le \frac{C_s}{\nu}\int_0^\infty x_4(t)^4\dd t
 \le \frac{C_s\eps^4}{\nu\mu}.\label{mean-max}
\end{equation}

Choose $A_s>1$ large enough that the exit bound is at most $\zeta/4$
when $R=A_s\zeta^{-1/2}\eps$. With $r_s$ already fixed, take the
smallness constant $c_s\le r_s/A_s$. Then
$\eps\le c_s\sqrt\zeta\sqrt{\nu\mu}$ ensures
$R\le r_s\sqrt{\nu\mu}$, as required for the moment estimates.
Use failure probability $\zeta/4$ for the exit estimate and for each of the two
Markov bounds from \eqref{simple-max} and \eqref{mean-max}.
All three bounds use the total initial norm, including the initial mean.
On the complement of these three exceptional events,
$\tau_R=\infty$ and
\eqref{zero-wave}--\eqref{zero-mean} hold, with constants of order
$C_s\zeta^{-1/2}$. The free heat contribution has supremum at most
$\|m_{\rm in}\|_{H^s}$ and satisfies
$\nu\int_0^\infty\|\partial_Ye^{\nu t\partial_{YY}}m_{\rm in}\|_{H^s}^2\dd t
\le\|m_{\rm in}\|_{H^s}^2/2$.
Together with \eqref{zero-small}, this proves \eqref{zero-total-mean}.

For the temperature phase, we first obtain
\eqref{zero-chart-reference} from
\eqref{zero-wave} and $\sup_{t\ge0}te^{-a\mu t}\le C_a\mu^{-1}$.
For the physical gradient, Brownian scaling and the reflection principle
on unit intervals give
\[
 \E\left[\left(\sup_{t\ge0}e^{-a\mu t/2}|W_t|\right)^2\right]
        \le C_a\mu^{-1}.
\]
We apply Markov's inequality with the remaining quarter of the failure probability.
Outside the resulting exceptional event, Sobolev embedding for $s\ge3$ and
$\nabla_t=(\partial_Y-\rho W_t\partial_Z,\partial_Z)$ yield
\[
 \|\nabla_t\chi(t)\|_\infty
 \le C_{s,\beta}t(1+\rho|W_t|)S(t)
 \le C_{s,\beta,\zeta}\eps\mu^{-1}
                              (1+\rho\mu^{-1/2}).
\]
Intersecting the resulting event with the preceding events proves \eqref{zero-chart-physical}.

It remains to prove decay of the mean and of the physical perturbation.
On the no-exit event, the weighted maximal bound is finite and the mean
flux belongs to $L^2_tH^s$ almost surely. For $t>T$, split the mean into
the heat evolution from $T$ and the forcing tail, bounded by
$C\nu^{-1/2}\|\mathcal F_{\rm m}\|_{L^2([T,\infty);H^s)}\to0$.
For each fixed splitting time, the heat evolution tends strongly to zero in $H^s(\R)$.
Taking first the large-time limit and then sending the splitting time to infinity gives $M(t)\to0$.
Since $|W_t|$ grows at most polynomially almost surely, the comparison
$\|f(t)\|_{H^s}\le C_s(1+\rho|W_t|)^s\|f(t)\|_{H^s_{\ad(t)}}$
implies decay of the simple part in physical $H^s$. The pullback
preserves means, so the complete perturbation decays.
\end{proof}

\subsection{Scalar mean normal forms}

For the scalar means, we refine \eqref{zero-mean} by cancelling their quadratic fluxes.
We take products in $(V,\Theta,Q)$,
whose $z$ averages agree with the physical averages. The derivatives are taken with
$\nabla_t=(\partial_Y-\rho W_t\partial_Z,\partial_Z)$.

\Needspace{5\baselineskip}
\begin{lemma}\label{lem-scalar-normal-form}
For the streamwise-independent solution, the following identities hold
almost surely in distributions on every finite time interval:
\begin{align}
 (\partial_t-\nu\partial_{YY})
 \left(\overline\Theta+
       \frac{1}{2\beta}\partial_Y\overline{\Theta^2}\right)
 &=-\frac{1}{2\beta}\partial_{YY}
           \overline{V^y\Theta^2}
   -\frac{\nu}{\beta}\partial_Y
           \overline{|\nabla_t\Theta|^2},
 \label{eq:theta-mean-normal-form}\\
 (\partial_t-\nu\partial_{YY})
 \left(\overline{Q}+
       \frac{1}{\beta}\partial_Y\overline{\Theta Q}\right)
 &=-\frac{1}{\beta}\partial_{YY}
           \overline{V^y\Theta Q}
   -\frac{2\nu}{\beta}\partial_Y
           \overline{\nabla_t\Theta\cdot\nabla_t Q}.
 \label{eq:q-mean-normal-form}
\end{align}
The right-hand sides contain cubic transport and viscosity-weighted
quadratic terms, with no stochastic integral or additional It\^o term.
\end{lemma}

\begin{proof}
The pullback removes the common scalar transport. With
$\Delta_t=\nabla_t\cdot\nabla_t$, we obtain
\begin{align*}
 \partial_t\Theta+V\cdot\nabla_t\Theta-\beta V^y
       &=\nu\Delta_t\Theta,\\*
 \partial_t Q+V\cdot\nabla_t Q
       &=\nu\Delta_t Q.
\end{align*}
Since $\nabla_t$ has constant coefficients at each time and
$\nabla_t\cdot V=0$, averaging a $\nabla_t$ divergence leaves only
the $\partial_Y$ derivative of its $y$ component. In particular,
$\overline{\Delta_t f}=\partial_{YY}\overline f$.
The product rule, with $L_\nu=\partial_t-\nu\partial_{YY}$, gives
\begin{align*}
 L_\nu\overline{\Theta^2}
 &=2\beta\overline{V^y\Theta}
   -\partial_Y\overline{V^y\Theta^2}
   -2\nu\overline{|\nabla_t\Theta|^2},\\
 L_\nu\overline{\Theta Q}
 &=\beta\overline{V^y Q}
   -\partial_Y\overline{V^y\Theta Q}
   -2\nu\overline{\nabla_t\Theta\cdot\nabla_t Q}.
\end{align*}
Since $\overline{V^y}=0$, the scalar mean equations are
$L_\nu\overline\Theta=-\partial_Y\overline{V^y\Theta}$ and
$L_\nu\overline{Q}=-\partial_Y\overline{V^y Q}$.
Combining the product identities with the mean equations gives the normal forms. The group approximation
of Proposition~\ref{prop:2d-global}, or the pulled-back product rule,
extends them to Sobolev solutions. Since Brownian paths are bounded on finite
intervals, the products and dissipative terms converge in distributions.
\end{proof}

\begin{corollary}
\label{cor-scalar-mean-uniform}
Let $s\ge3$ be an integer, $\nu,\beta>0$, $\rho\ge0$, and
$\zeta\in(0,1)$.  Suppose the deterministic initial state for
\eqref{eq:2d-subsystem} belongs to $H^{s+2}(\mathcal D)$, its velocity
is divergence-free, and all its $z$ means are zero. Put
$\eps=\|(v_{\rm in},\theta_{\rm in},q_{\rm in})\|_{H^{s+2}}$ and
$\mu=\nu+\rho\sqrt\nu$.
There are $c_{s,\zeta}>0$ and $C_{s,\beta,\zeta}>0$, independent of
$\nu$ and $\rho$, such that, if $\eps\le c_{s,\zeta}\sqrt{\nu\mu}$,
then, with probability at least $1-\zeta$,
\begin{equation}\label{eq:scalar-mean-uniform}
 \sup_{t\ge0}
 \left(\|\overline\theta(t)\|_{H^s(\R)}
       +\|\overline q(t)\|_{H^s(\R)}\right)
 \le C_{s,\beta,\zeta}\eps^2.
\end{equation}
This holds on the event of Proposition~\ref{thm-zero} at order $s+2$.
\end{corollary}

\begin{proof}
Set $j=s+2$ and $S_j=\|\widetilde{Y^\sh}\|_{H^j}$.
On the event of Proposition~\ref{thm-zero} at order $j$, the bounds and smallness
condition give, with $K_{s,\zeta}$ independent of $\nu,\rho$,
\begin{align*}
 A_j:=\sup_{t\ge0}\|Y^\sh(t)\|_{H^j}
 &\le K_{s,\zeta}\eps,\\
 \int_0^\infty S_j(t)^2\dd t
 &\le \frac{K_{s,\zeta}\eps^2}{\mu},\\
 \nu\int_0^\infty
       \|\nabla_t Y^\sh(t)\|_{H^{s+1}}^2\dd t
 &\le K_{s,\zeta}\eps^2.
\end{align*}
Indeed, the gradients of the simple and mean components are orthogonal in the reference norm.
Squaring \eqref{zero-wave} and \eqref{zero-mean} at order $j$ gives
\[
 \nu\int_0^\infty\|\nabla_t Y^\sh\|_{H^j}^2\dd t
 \le K_{s,\zeta}
       \left(\eps^2+\frac{\eps^4}{\nu\mu}\right)
 \le K_{s,\zeta}\eps^2,
\]
after increasing $K_{s,\zeta}$; for means, $\nabla_t=(\partial_Y,0)$.

We first bound the quadratic corrections and sources using the Sobolev algebra
property and the boundedness of normalized averaging:
\begin{align*}
 \|\partial_Y\overline{\Theta^2}\|_{H^s}
 +\|\partial_Y\overline{\Theta Q}\|_{H^s}
 &\le C_s A_j^2,\\
 \|\partial_Y\overline{|\nabla_t\Theta|^2}\|_{H^s}
 +\|\partial_Y\overline{
            \nabla_t\Theta\cdot\nabla_t Q}\|_{H^s}
 &\le C_s\|\nabla_t Y^\sh\|_{H^{s+1}}^2.
\end{align*}
For the cubic terms, $V^y=\widetilde V^y$. Hence, in each flux, the term containing two
scalar means has zero average. Expanding the remaining terms gives
\[
 \|\overline{V^y\Theta^2}\|_{H^{s+2}(\R)}
 +\|\overline{V^y\Theta Q}\|_{H^{s+2}(\R)}
 \le C_s A_j S_j^2.
\]
The cubic terms in
\eqref{eq:theta-mean-normal-form}--\eqref{eq:q-mean-normal-form}
have $L^1_tH^s$ norm at most
$C_{s,\beta,\zeta}\eps^3/\mu$, and the viscosity-weighted
quadratic terms have $L^1_tH^s$ norm at most
$C_{s,\beta,\zeta}\eps^2$.

We apply Duhamel's formula and heat contraction on $H^s(\R)$ to both
normal forms. The initial means vanish, and the quadratic corrections
at the initial time and at time $t$ are bounded by $C_{s,\beta,\zeta}\eps^2$ in
$H^s$. Combining these bounds with the source estimates yields
\[
 \sup_{t\ge0}
 \left(\|\overline\Theta(t)\|_{H^s}
       +\|\overline{Q}(t)\|_{H^s}\right)
 \le C_{s,\beta,\zeta}
       \left(\eps^2+\frac{\eps^3}{\mu}\right).
\]
Since $\mu\ge\nu$, the smallness condition gives
$\eps/\mu\le c_{s,\zeta}\sqrt{\nu/\mu}\le c_{s,\zeta}$.
The cubic contribution therefore has the same quadratic bound. Since the pullback
preserves scalar means, this proves
\eqref{eq:scalar-mean-uniform}.
\end{proof}

\begin{remark*}
With zero initial means and two additional Sobolev derivatives in the initial data,
the scalar mean supremum improves from $O(\eps^2/\sqrt{\nu\mu})$ to
$O(\eps^2)$ under the smallness condition of Proposition~\ref{thm-zero}.
The constant in \eqref{eq:scalar-mean-uniform} may depend on $\beta$ and $\zeta$.
\end{remark*}

\Needspace{5\baselineskip}
\section[Transverse perturbations and three-dimensional stability]{Transverse perturbations and 
three-dimensional stability}\label{sec:3d-stability}

We now use the estimates for the streamwise-independent system to
control its interaction with the fast modes. We first consider transverse
linear perturbations, and then close the estimates for the full
nonlinear system.

\subsection{Nonlinear product estimates}

We work in the common affine coordinates $X=x-ty$, $Y=y$, $Z=z-\rho W_ty$,
with $\nabla_t=(\partial_X,\partial_Y-t\partial_X-\rho W_t\partial_Z,
                         \partial_Z)$.
Write $U=(u,\theta)$, with componentwise pullback $U^\sh$,
and denote the reference $H^s$ norm by $\norm{\cdot}_s$.
The symbol of $\Lambda_t=|\nabla_t|$ is
$r_t(\xi)=|(k,\eta-kt-\rho\ell W_t,\ell)|$, where $\xi=(k,\eta,\ell)$.
For $k\ne0$, set $E=\mathcal E_s$, $D=\mathcal D_s$ and
$L=\mathcal L_s$ as in \eqref{eq:lp-fast-field-energy}; then
$D\ge E$ and $L\ge\nu^{1/3}E$. For the full $x$-independent field, we set
$Y_0=(v_0,\theta_0,q_0)$, $v_0=(u_0^2,u_0^3)$,
$q_0=u_0^1+\beta^{-1}\theta_0$,
$E_0=\norm{Y_0^\sh}_s^2$, and $G_0=\norm{\nabla_tY_0^\sh}_s^2$.
Thus $E_0$ includes both zero sectors.
Constants may depend on $s,\beta$.

\begin{lemma}\label{lem:joint-products}
For $s>3/2$, the physical half derivative satisfies
\begin{equation}\label{eq:joint-half-product}
 \norm{\Lambda_t^{1/2}(AB)}_s
 \le C_s\left(\norm A_s\norm{\Lambda_t^{1/2}B}_s
                  +\norm{\Lambda_t^{1/2}A}_s\norm B_s\right),
\end{equation}
uniformly in $t$ and $W_t$. For a regular solution of
\eqref{eq:model}, let $\mathcal N_0$ and $\mathcal N_f$ be the nonlinear
contributions to $\dd E_0$ and $\dd E$. If $s\ge3$ and $\beta>1/2$, then
\begin{equation}\label{eq:joint-nonlinearity}
 |\mathcal N_0|+|\mathcal N_f|
       \le C_{s,\beta}\sqrt{E_0+E}\,(G_0+D).
\end{equation}
\end{lemma}

\begin{proof}
The affine frequency map gives
$r_t(\xi+\xi')^{1/2}\le r_t(\xi)^{1/2}+r_t(\xi')^{1/2}$.
Combining this inequality with $\av{\xi+\xi'}^s\le C_s(\av\xi^s+\av{\xi'}^s)$,
Young's convolution inequality and
$\sum_{k,\ell}\int_{\R}\av{k,\eta,\ell}^{-2s}\dd\eta<\infty$
yields \eqref{eq:joint-half-product}. Omitting the factor
$r_t^{1/2}$ in the convolution argument gives the uniform algebra bound
in the reference $H^s$ norm.

Set $n=\sqrt E$, $g=\sqrt D$, $a=\sqrt{E_0}$ and $b=\sqrt{G_0}$.
We reconstruct $U_f^\sh$ from $(v,\theta)$ and
$p=b_{\rm pol}+\ell\theta/(\beta h)$. The polarization,
reconstruction and Leray multipliers are bounded, so coercivity gives
\begin{equation*}
 \norm{U_f^\sh}_s\le C_\beta n,
  \quad
 \norm{\Lambda_tU_f^\sh}_s\le C_\beta g,
  \quad g\ge n.
\end{equation*}
The half-weighted wave norm is at most $C_\beta n$; interpolation gives
\[
 \norm{\Lambda_t^{1/2}p}_s
       \le\norm p_s^{1/2}\norm{\Lambda_tp}_s^{1/2}
       \le C_\beta\sqrt{ng}.
\]
Since $g\ge n$, reconstruction yields
$\norm{\Lambda_t^{1/2}U_f^\sh}_s\le C_\beta\sqrt{ng}$,
$\norm{U_0^\sh}_s\le C_\beta a$, and
$\norm{\Lambda_t^{1/2}U_0^\sh}_s\le C_\beta\sqrt{ab}$.
Since this interpolation is homogeneous, it also applies to the double
mean without Poincar\'e's inequality.

We keep the divergence outside the flux and bound the fast wave-energy pairing
by $C_\beta g\norm{\Lambda_t^{1/2}\mathcal H_f}_s$.
In this pairing, the energy weight and the divergence supply $r_t^2$, which we split as $r_t^{3/2}$
on the wave factor and $r_t^{1/2}$ on the flux. For the passive pairing,
we need $g\norm{\mathcal H_f}_s$, which satisfies the same bound since
$r_t\ge1$ on $k\ne0$. Keeping the divergence outside the flux thus avoids
commuting multipliers through products.

By \eqref{eq:joint-half-product}, the fast--fast and zero--fast fluxes
have half-weighted norms bounded by $Cn\sqrt{ng}$ and
$C[a\sqrt{ng}+n\sqrt{ab}]$, respectively. Hence
\[
 |\mathcal N_f|
 \le C_{s,\beta}\left[n^{3/2}g^{3/2}
             +a\sqrt n\,g^{3/2}+\sqrt a\,ng\sqrt b\right].
\]
Let $R=\sqrt{E_0+E}$. Since $a,n\le R$ and $g\ge n$, the first two terms in this bound
are bounded by $CRg^2$. Using $n^2\le Rg$, the last is bounded by
$CRg^{3/2}b^{1/2}\le CR(g^2+b^2)$.

We next estimate the zero-output flux, which contains only zero--zero and fast--fast products.
Set $M=U_{00}$ and $S=\norm{U_0^\sh-M}_s$. Incompressibility and
finite energy give $M_u^2=0$. The mean--mean fluxes therefore have zero
divergence, and the remaining zero--zero flux is bounded by $CaS$, where
$S\le C_\beta b$ since $\ell\ne0$. The fast--fast flux is bounded by
$Cn^2$. In the $q_0$ equation, $x$ averaging removes the pressure
derivative, and the shear term cancels the temperature source. Pairing the
resulting divergence flux with $\nabla_tY_0^\sh$ gives
$|\mathcal N_0|\le C_{s,\beta}b(ab+n^2)
             \le C_{s,\beta}R(b^2+g^2)$,
using $n^2\le Rg$. Combining the fast and zero-mode bounds proves
\eqref{eq:joint-nonlinearity}.
\end{proof}

\subsection{Separated mixed interactions}

To retain the separate sector sizes, we refine the joint bound and
estimate the fast self-interactions and zero--fast interactions separately.

Put $\lambda=\nu^{1/3}$, $\mu=\nu+\rho\sqrt\nu$,
$A_0=\norm{U_0^\sh}_{H^s}$, and
$H_0=\norm{\Lambda_t^{1/2}U_0^\sh}_{H^s}$.
Homogeneous interpolation and reconstruction give
\begin{equation}\label{eq:separated-zero-half}
 A_0^2\le C_\beta E_0,
  \quad H_0^2\le C_\beta\sqrt{E_0}\sqrt{G_0}.
\end{equation}
These bounds include the double mean and do not use a spectral gap.

We introduce a wave weight $1\le n\le e^\pi$ and replace $F+F_p$ by
$nF+F_p$ in $E,D,L$. The weight acts at the output frequency.
We decompose the instantaneous Fourier coefficients as
$\widehat{U_P^\sh}=(e_vp,0)$ and
$\widehat{U_W^\sh}=(e_gv-\ell\theta e_v/(\beta h),\theta)$.
We evaluate the output factors at $\eta_t$ and set
$B_W=\|U_W^\sh\|_{H^s}^2$ and
\[ Z_W=\sum_{k\ne0,\ell}\int_\R\av\xi^{2s}\frac{|k|h}{r_t^2}nF\dd\eta. \]

\begin{lemma}\label{lem:separated-mixed}
Let $s\ge3$ and $\beta>1/2$. The fast self-interaction and mixed
zero--fast contributions satisfy
\begin{align}
 |\mathcal N_{ff}|
 &\le C_{s,\beta}(ED^{1/2}+B_W^{1/2}E^{1/4}D^{3/4}),\label{eq:separated-self}\\
 |\mathcal N_{0f}|
 &\le C_{s,\beta}(A_0E^{1/4}D^{3/4}
                         +H_0E^{1/2}D^{1/2}).\label{eq:separated-mixed}
\end{align}
Consequently, for every $\eta>0$,
\begin{align}
 |\mathcal N_{ff}|
 &\le\eta\nu D+C_{s,\beta,\eta}
                              (\nu^{-1}E^2+\nu^{-3}E^2Z_W),\label{eq:separated-self-young}\\
 |\mathcal N_{0f}|
 &\le\eta\nu D+C_{s,\beta,\eta}
                         (\nu^{-3}A_0^4+\nu^{-1}H_0^2)E.
                         \label{eq:separated-mixed-young}
\end{align}
The constants are uniform in $t$ and $W_t$.
\end{lemma}

\begin{proof}
The half-weighted wave norm is bounded by $C_\beta\sqrt E$, while the ordinary
wave norm is bounded by $\sqrt{B_W}$. Applying \eqref{eq:joint-half-product}, we bound the wave--wave
and wave--passive fluxes by
$C(E+B_W^{1/2}E^{1/4}D^{1/4})$.
It remains to estimate two passive inputs with reference frequencies $\xi_a,\xi_b$
and a fast output $\xi_c=\xi_a+\xi_b$. We write
$K_j=(k_j,\ell_j)$, $h_j=|K_j|$, and
$\Delta=k_a\ell_b-\ell_ak_b$. The affine map preserves frequency addition.
With the convention $\partial_j\leftrightarrow
i\xi_j$, the ordered velocity source is
$-i(e_{v,a}\cdot\xi_{b,t})\mathcal P_c e_{v,b}\,p_ap_b$, where
$\mathcal P_c=\mathcal P(\xi_{c,t})$ and
$e_{v,a}\cdot\xi_{b,t}=-\frac{\Delta}{h_a}$.
Since $e_{g,c}$ belongs to the range of $\mathcal P_c$,
$e_{g,c}\cdot \mathcal P_c e_{v,b}=-\frac{\eta_{c,t}\Delta}{r_ch_ch_b}$.
It follows that the wave source is
$-i\frac{\eta_{c,t}\Delta^2}{r_ch_ch_ah_b}\,p_ap_b$.
Symmetrizing in $a,b$ leaves this coefficient unchanged and gives the
passive output
$\frac{i\Delta(h_b^2-h_a^2)}{2h_ch_ah_b}\,p_ap_b$.
The passive--passive thermal source vanishes, so this is also the
contribution to the $p$ equation.

For the wave coefficient, put $\Delta=h_ah_b\sin\phi$. The identities
$K_c=K_a+K_b$ and
$h_c\ge h_b|\sin\phi|$, $h_c\ge h_a|\sin\phi|$ give
\[
 \frac{|\eta_{c,t}|\Delta^2}{r_ch_ch_ah_b}
 \le\frac{h_ah_b\sin^2\phi}{h_c}
 \le\min(h_a,h_b).
\]
The bound on the wave coefficient allows us, in the $H^s$ convolution
estimate, to place the derivative on the lower-weight factor in
$\av{\xi_c}^s\le C_s(\av{\xi_a}^s+\av{\xi_b}^s)$.
The estimate
\[
 \|h\widehat p\|_{L^1}
 \le\left(\sum_{k,\ell}\int_\R h^2\av{k,\eta,\ell}^{-2s}\dd\eta\right)^{1/2}
       \|p\|_{H^s}
 \le C_s\|p\|_{H^s}, \quad s>5/2,
\]
and Young's convolution inequality bound the wave source in $H^s$ by
$CE$. The wave-energy pairing is therefore at most $CED^{1/2}$, since
division by the wave dissipation factor leaves $r_c^{-1/2}\le1$.
For the passive output, the flux $u_P\otimes u_P$ has $H^s$ norm at most
$C\|p\|_{H^s}^2$. Keeping its divergence outside and pairing with the
output gradient gives $CED^{1/2}$. Bounded output multipliers then give
\eqref{eq:separated-self}.

We next relate the ordinary wave norm to the positive wave term.
Coercivity and Cauchy--Schwarz give
\[
 B_W\le C_\beta\sum_{k\ne0,\ell}\int_\R
                      \av\xi^{2s}r_t^{-1}nF\dd\eta,
  \quad
 B_W^2\le C_\beta E\sum_{k\ne0,\ell}\int_\R
                      \av\xi^{2s}r_t^{-2}nF\dd\eta
          \le C_\beta E Z_W,
\]
where the last inequality uses $|k|h\ge1$. Hence
\[
 B_W^{1/2}E^{1/4}D^{3/4}
 \le\eta\nu D+C_{\beta,\eta}\nu^{-3}B_W^2E
 \le\eta\nu D+C_{\beta,\eta}\nu^{-3}E^2Z_W.
\]
Combining the latter estimate with Young's inequality for $ED^{1/2}$ yields
\eqref{eq:separated-self-young}. Absorbing $Z_W$ requires
$E^2\lesssim\nu^3$, which gives the scale $\alpha\lesssim\nu^{3/4}$.
For the zero--fast flux, Lemma~\ref{lem:joint-products} gives
$C[A_0(E^{1/2}+E^{1/4}D^{1/4})+H_0E^{1/2}]$.
Multiplying by $D^{1/2}$ and using $D\ge E$, we obtain
\eqref{eq:separated-mixed}; Young's inequality then gives
\eqref{eq:separated-mixed-young}.
\end{proof}

We next construct a bounded wave weight that supplies the positive
term needed to absorb the fast self-interaction.

\begin{lemma}\label{lem:critical-poisson}
For $k\ne0$, set $d=\rho^2\ell^2/(2|k|)$ and define
\[
 H(\eta)=\int_0^\infty e^{-z}
       \frac{h}{h^2+(\eta+\operatorname{sgn}(k)dz)^2}\dd z,
  \quad
 g(\eta)=
 \begin{cases}\displaystyle\int_{-\infty}^{\eta}H(y)\dd y,&k>0,\\
              \displaystyle\int_{\eta}^{\infty}H(y)\dd y,&k<0.
 \end{cases}
\]
Then $0\le g\le\pi$, $|g'|\le h^{-1}$, $|g''|\le Ch^{-2}$, and
\begin{equation}\label{eq:critical-poisson}
 \frac{\rho^2\ell^2}{2}g''-kg'=-\frac{|k|h}{r^2}.
\end{equation}
Define $E_g,D_g,L_g$ by replacing $F$ with $e^gF$
in \eqref{eq:lp-fast-field-energy}, keeping $F_p$ unchanged.
Let $Z_g=Z_W$ for $n=e^g$. The linear contribution satisfies
\begin{equation}\label{eq:critical-linear-drift}
 \dd E_g+\nu D_g\dd t+L_g\dd t+Z_g\dd t
 \le C_\beta\rho^2E_g\dd t+\rho G_g\dd W_t,
  \quad |G_g|\le C_\beta E_g.
\end{equation}
\end{lemma}
\begin{proof}
By definition, $g'=\operatorname{sgn}(k)H$; integration by parts gives
$\operatorname{sgn}(k)dH'=H-h/r^2$.
Thus \eqref{eq:critical-poisson} also holds when $\rho\ell=0$.
The integral formula also gives the stated bounds.
For $n=e^g$, the additional It\^o drift in $\dd(nF)$ is
\[
 -\frac{|k|h}{r^2}nF+
 \frac{\rho^2\ell^2}{2}(g')^2nF+
 \rho^2\ell^2ng'X^*(\mathsf N'+\eta\mathsf N/r^2)X.
\]
The last two terms are bounded by $C_\beta\rho^2nF$, since
$|\ell|\le h$, $|X^*\mathsf N'X|\le C_\beta h^{-1}F$, and
$|\eta|/r^2\le(2h)^{-1}$. The additional martingale is
$-\rho\ell g'nF\dd W_t$. Applying Lemma~\ref{lem:lp-fast-Ito} and
integrating over Fourier modes gives \eqref{eq:critical-linear-drift}.
The truncation argument of Proposition~\ref{prop:lp-fast-field}
extends this estimate to finite-energy fields.
\end{proof}

\subsection{A compensated fast-mode estimate}

The mixed estimate contains a coefficient whose time integral is
controlled by the zero-mode energy and dissipation. To retain fast-mode
decay in the presence of this coefficient, we incorporate its accumulated
contribution into an adapted weight.

\begin{lemma}\label{lem:compensated-fast}
Let $E,D,L$ be nonnegative adapted processes, with $E$ a continuous
semimartingale and $L\ge\lambda E$. Suppose that, on $[0,\infty)$,
\begin{equation}\label{eq:compensated-fast-drift}
 \begin{aligned}
 \dd E+\frac\nu2D\dd t+\frac14L\dd t+\frac\lambda2E\dd t
       &\le a_tE\dd t+\rho G_t\dd W_t,\\
 |G_t|&\le C_\beta E,
 \end{aligned}
\end{equation}
where $a_t\ge0$ is progressively measurable and locally integrable, and
$G$ is predictable. Let $\alpha^2=E(0)$ be deterministic and assume
\begin{equation}\label{eq:compensated-fast-accumulator}
 C_t=\int_0^t a_r\dd r\le\frac\lambda8t+\Gamma
\end{equation}
for a deterministic $\Gamma\ge0$. If
$\rho^2\le c_\beta\lambda$ with a sufficiently small constant,
then for $q=1,2,4$,
\begin{equation}\label{eq:compensated-fast-moments}
 (\E E(t)^q)^{1/q}\le e^\Gamma\alpha^2e^{-\lambda t/8},
\end{equation}
and
\begin{align}
 \E\left[\left(\sup_{t\ge0}e^{\lambda t/8}E(t)\right)^2\right]
          &\le2e^{2\Gamma}\alpha^4,\label{eq:compensated-fast-maximal}\\
 \E\int_0^\infty e^{\lambda t/8}(\nu D+L)\dd t
          &\le4e^\Gamma\alpha^2.\label{eq:compensated-fast-dissipation}
\end{align}
For a stopping time $\tau$, the moment conclusion is
\[
 \E\left[e^{q\lambda(t\wedge\tau)/8}
                   E(t\wedge\tau)^q\right]
        \le e^{q\Gamma}\alpha^{2q}.
\]
The maximal estimate uses the weight at $t\wedge\tau$, with dissipation integrated
over $[0,\tau]$. If the homogeneous linear continuation satisfies
\eqref{eq:compensated-fast-drift} after exit, the bounds with time $t$
hold for all $t\ge0$. No independence of $C$ and $W$ is required.
\end{lemma}

\begin{proof}
Set $V=e^{-C}E$. Since $e^{-C}$
is adapted and has finite variation,
\eqref{eq:compensated-fast-drift} gives
\[
 \dd V+\frac\nu2e^{-C}D\dd t
          +\frac14e^{-C}L\dd t
          +\frac\lambda2V\dd t
       \le\rho\widetilde G\dd W_t,
  \quad |\widetilde G|\le C_\beta V.
\]
We choose the noise constant small enough to control the fourth moment.
By It\^o's formula, $e^{q\lambda t/4}V^q$ is a nonnegative local
supermartingale for $q=1,2,4$: the correction
$C_\beta q(q-1)\rho^2V^q$ is absorbed by $q\lambda V^q/4$.
Localization and Fatou's lemma yield
$\E V(t)^q\le\alpha^{2q}e^{-q\lambda t/4}$.
Using \eqref{eq:compensated-fast-accumulator}, we obtain
\eqref{eq:compensated-fast-moments}.

We use the maximal inequality for the fourth-power supermartingale.
Integrating the resulting tail bound to control the second moment gives
$\E\left[\left(\sup_{t\ge0}e^{\lambda t/4}V(t)\right)^2\right]
                                     \le2\alpha^4$.
Since $e^{\lambda t/8}E(t)\le e^\Gamma e^{\lambda t/4}V(t)$,
this proves \eqref{eq:compensated-fast-maximal}.
For the dissipation estimate, we multiply the first-moment inequality by
$e^{\lambda t/4}$ and obtain
$\E\int_0^\infty e^{\lambda t/4-C_t}
                    (\nu D+L)\dd t\le4\alpha^2$.
Since the resulting weight is at least $e^{-\Gamma}e^{\lambda t/8}$,
\eqref{eq:compensated-fast-dissipation} follows by localization.
Replacing $t$ by $t\wedge\tau$ gives the stopped bounds. The linear
continuation has $a_t=0$ after exit, so the argument holds for all $t\ge0$.
\end{proof}

\subsection{Transverse perturbations of random zero-mode solutions}

We now fix a random streamwise-independent background and consider its
transverse linearization. Throughout this subsection, we write $E,D,L$ for the fast energies of the transverse linearization in
\eqref{eq:lp-fast-field-energy}, and $E_0,G_0$ for the background energy
and dissipation. All products are formed using the componentwise pullbacks.

\begin{proposition}\label{prop:transverse-stability}
Fix an integer $s\ge3$, $\beta>1/2$, $0<\nu\le1$ and
$\zeta\in(0,1)$, and put $\lambda=\nu^{1/3}$.
Let a global nonlinear streamwise-independent solution have zero
initial double mean. Suppose, on an event of probability at least
$1-\zeta/2$, that
$\sup_{t\ge0}E_0(t)+\nu\int_0^\infty G_0(t)dt
                    \le K\varepsilon_b^2$, with $K\ge1$.
If $\varepsilon_b\le c_{s,\beta,K}\nu^{5/6}$ and
$\rho^2\le c_{s,\beta}\lambda$, its transverse linearization with
deterministic initial fast energy $E(0)=\alpha^2$ satisfies, with
probability at least $1-\zeta$,
\begin{equation}\label{eq:transverse-stability-decay}
 \sup_{t\ge0}e^{\lambda t/8}E(t)
 +\int_0^\infty e^{\lambda t/8}(\nu D+L)dt
                         \le C_{\beta,\zeta}\alpha^2.
\end{equation}
The fast size $\alpha$ is independent of $\varepsilon_b$.
The estimate uses only the stated background $H^s$ and dissipation bounds.
Proposition~\ref{thm-zero} provides such solutions under
$\varepsilon_b\lesssim_{s,\beta,\zeta}
\min\{\nu^{5/6},\sqrt{\nu\mu}\}$, where
$\mu=\nu+\rho\sqrt\nu$.
\end{proposition}

\Needspace{7\baselineskip}
\begin{proof}
For $\varepsilon_b=0$, the homogeneous fast estimate applies.
Otherwise, we stop the background~at
\[
 \sigma=\inf\left\{t:E_0(t)+\nu\int_0^tG_0(r)dr
                               \ge2K\varepsilon_b^2\right\}.
\]
On the event where the background bound holds, $\sigma=\infty$.
To estimate the mixed contribution, we use reconstruction and homogeneous
interpolation to obtain $A_0^2\le C_\beta E_0$
and $H_0^2\le C_\beta\sqrt{E_0G_0}$. Keeping the two terms of
Lemma~\ref{lem:joint-products} separate and placing the output derivative
on the fast dissipation factor, we obtain
\[
 \begin{aligned}
 |\mathcal N_{0f}|&\le C_{s,\beta}D^{1/2}
       [(A_0+H_0)E^{1/2}+A_0E^{1/4}D^{1/4}]\\
 &\le\frac\nu2D
       +C_{s,\beta}[\nu^{-1}(A_0^2+H_0^2)+\nu^{-3}A_0^4]E.
 \end{aligned}
\]
Before $\sigma$, the terms involving the background Sobolev norm are bounded by
$C_{s,\beta}[K\varepsilon_b^2/\nu+K^2\varepsilon_b^4/\nu^3]$.
By the smallness assumptions, we can absorb these terms together with
the linear noise drift and obtain
\[
 \begin{gathered}
 dE+\frac\nu2Ddt+\frac14Ldt+\frac\lambda2Edt
       \le b(t)Edt+\rho G_tdW_t,\\
 b=C_{s,\beta}\nu^{-1}H_0^2,
 \quad |G_t|\le C_\beta E.
 \end{gathered}
\]
For $u=t\wedge\sigma$, Cauchy--Schwarz and the stopped background energy bound
give
\[
 A_u:=\int_0^ub(r)dr
 \le C_{s,\beta}K\varepsilon_b^2\nu^{-3/2}\sqrt u
 \le\frac\lambda8u+\Gamma,
  \quad
 \Gamma=C_{s,\beta}K^2\varepsilon_b^4/(\nu^3\lambda)\le1.
\]
We apply Lemma~\ref{lem:compensated-fast} with weights stopped at
$u=t\wedge\sigma$. For $V_t=e^{\lambda u/4-A_u}E(u)$, the fourth power
is a nonnegative supermartingale after localization. The maximal and
first-moment estimates give
\[
 \mathbb E\bigl[\sup_t V_t^2\bigr]\le2\alpha^4,
  \quad
 \mathbb E\int_0^\sigma e^{\lambda t/4-A_t}
                        (\nu D+L)dt\le4\alpha^2.
\]
Since $A_u\le\lambda u/8+1$, these estimates control the stopped quantities in
\eqref{eq:transverse-stability-decay}. We apply Markov's inequality with failure
probability $\zeta/2$ and intersect the resulting event with the
background event to obtain the claimed bound.

It remains to verify that the zero-mode result supplies the required
background bound. We remove the mean--mean flux from the zero-mode
energy identity. The remaining flux is bounded by $C_\beta\sqrt{E_0}\,S$,
so
$dE_0+\nu G_0dt\le C_{s,\beta}\nu^{-1}E_0S^2dt$.
On the event in Proposition~\ref{thm-zero}, we have
$\sup E_0\le C_\zeta\varepsilon_b^2$ and
$\int S^2dt\le C_\zeta\varepsilon_b^2/\mu$.
Consequently
$\nu\int G_0dt\le C_{s,\beta,\zeta}
[\varepsilon_b^2+\varepsilon_b^4/(\nu\mu)]$,
as required when
$\varepsilon_b\lesssim_{s,\beta,\zeta}\sqrt{\nu\mu}$.
\end{proof}

\subsection{Full nonlinear stability}

In the full system, the fast modes generate a zero field, which in turn
feeds back into the fast equations. We control this coupling by closing
both energy bounds at a common stopping time.

For every local strong solution of \eqref{eq:model}, averaging in
$x,z$ gives $(\partial_t-\nu\partial_{yy})f_{00}
 =-\partial_y(u^2f)_{00}$ for $f=u^1,u^3,\theta$.
Here $u^2_{00}=0$. Averaging removes both affine transports, the buoyancy
term and the lift-up term, while the averaged Leray projection preserves the first
and third velocity components.

\begin{lemma}\label{lem:fast-source-convolution}
Let $0<\mu\le\lambda$, $c>0$, and
$0<b\le\min\{c/2,1/16\}$. Then
\begin{equation*}
 \sup_{t\ge0}e^{b\mu t}\int_0^t
 \frac{e^{-c\mu(t-r)}e^{-\lambda r/8}}{\sqrt{\nu(t-r)}}\dd r
 \le\frac{C}{\sqrt{\nu\lambda}}.
\end{equation*}
\end{lemma}

\begin{proof}
The weighted integrand is bounded by
$[\nu(t-r)]^{-1/2}e^{-\lambda r/16}$.
With $z=\lambda t$, $u=\lambda r$, the integral of this upper bound becomes
$(\nu\lambda)^{-1/2}\int_0^z(z-u)^{-1/2}e^{-u/16}\dd u$.
For $z\le1$, the rescaled integral is at most $2\sqrt z$. For $z\ge1$, we split
the rescaled integral at $z/2$ and bound it by $16\sqrt{2/z}+2\sqrt z\,e^{-z/32}$.
This bound is uniform in $z$.
\end{proof}

\begin{lemma}\label{lem:physical-fast-damping}
Let $s\ge3/2$, and use the divergence-free fast polarization in
\eqref{eq:intro-polarization}--\eqref{eq:intro-fast-energies}.
With $B_t=1+\rho|W_t|$, one has
\begin{equation}\label{eq:physical-fast-frequency-bound}
 \begin{aligned}
 \norm{u_{\mathrm{nz}}^2(t)}_{L^2}^2&\le C\av{t}^{-3}B_t^3 E_{\mathrm{nz}}(t),\\
 \norm{\theta_{\mathrm{nz}}(t)}_{L^2}^2&\le C\av{t}^{-1}B_t E_{\mathrm{nz}}(t).
 \end{aligned}
\end{equation}
If $\lambda>0$ and $\rho\le R\sqrt\lambda$, then
\begin{equation}\label{eq:physical-fast-brownian-event}
 \Prob\left(\sup_{t\ge0}e^{-\lambda t/16}B_t^3
              \le C_{R,\zeta}\right)\ge1-\zeta/2,
  \quad 0<\zeta<1,
\end{equation}
where $C_{R,\zeta}$ is independent of $\lambda,\rho$.
\end{lemma}

\begin{proof}
We write $d=\av\xi$ for $k\ne0$. Since $r_t\ge h\ge1$ and $d\ge h$,
the characteristic identity gives
$t\le t|k|\le|\eta|+|\eta_t|+\rho|\ell||W_t|\le r_t+B_td$.
Hence $\av t\le3B_tr_td$ and $\av t\le3B_t\max\{r_t,d\}$.
Combining these bounds with $h^2/(r_t^3d^3)\le\min\{r_t^{-3},d^{-3}\}$ gives
$\frac1{r_td}\le\frac{3B_t}{\av t}$ and
$\frac{h^2}{r_t^3d^3}\le\frac{27B_t^3}{\av t^3}$.
Since $e_v^2=0$ and $e_g^2=h/r_t$,
$\widehat{(u_{\mathrm{nz}}^2)^\sh}=(h/r_t)v$.
For $2s\ge3$, the squared velocity and temperature coefficients, divided
by the energy density, are bounded by $h^2/(r_t^3d^3)$ and $1/(r_td)$, respectively.
Parseval's identity and the $L^2$-isometry of the pullback give
\eqref{eq:physical-fast-frequency-bound}.

For the probability bound, Brownian scaling gives a standard Brownian motion
$\widetilde W_u=\sqrt\lambda W_{u/\lambda}$. The supremum in \eqref{eq:physical-fast-brownian-event}
is at most $Z_R=\sup_{u\ge0}e^{-u/16}(1+R|\widetilde W_u|)^3$.
The Brownian maximal inequality gives
\[
 \E Z_R\le\sum_{n\ge0}e^{-n/16}
 \E(1+R\sup_{u\le n+1}|\widetilde W_u|)^3
 \le C_R\sum_{n\ge0}e^{-n/16}[1+(n+1)^{3/2}]<\infty.
\]
Markov's inequality with $C_{R,\zeta}=2\E Z_R/\zeta$ proves the claim.
\end{proof}

\begin{proof}[Proof of Theorem~\ref{thm:full3d-stability}]
By Lemma~\ref{lem:mixed-local}, the solution of
Proposition~\ref{thm:affine-local-strong} has continuous polarization energies
and satisfies the nonlinear energy identities up to locally bounded
stopping times. We use the bounded correction of Lemma~\ref{lem:critical-poisson} and write
$E=E_g$, $D=D_g$, $L=L_g$, $Z=Z_g$. Coercivity gives
$E\asymp_\beta E_{\mathrm{nz}}$, $D\asymp_\beta D_{\mathrm{nz}}$, and
$L\asymp_\beta L_{\mathrm{nz}}$.
Hence $E(0)\le C_\beta\alpha^2$, and the estimates transfer to
\eqref{eq:intro-fast-energies} with fixed equivalence constants.

Set $\sigma_0=\min\{\nu^{5/6},\sqrt{\nu\mu}\}$.
The fast-mode data condition gives
$\frac{\alpha^2}{\sqrt{\nu\lambda}}
       \le c_{s,\beta}^2\zeta^{1/2}\sigma_0
       \le c_{s,\beta}^2\zeta^{1/4}\sigma_0$.
We shall choose $K\ge C_{s,\beta}/\sqrt\zeta$ at the end of the
stopping-time argument. For this choice, a sufficiently small
constant in \eqref{eq:full3d-small} gives
\begin{equation}\label{eq:split-stopping-smallness}
 \sqrt K\alpha\le c_*\nu^{3/4}, \quad
 \sqrt K\delta\le c_*\sigma_0,
\end{equation}
where $c_*>0$ is fixed and depends only on $s,\beta$.
The data conditions also give $\delta\le C\eps$.

We define $Q_0(t)=E_0(t)+\nu\int_0^tG_0(r)\dd r$ and
\[
 \tau=\inf\{t:Q_0(t)\ge K\delta^2\text{ or }E(t)\ge K\alpha^2\}.
\]
If $\alpha=0$, we omit the fast threshold: uniqueness and invariance of
$k=0$ imply that the fast field remains zero. If $\delta=0$, both fields vanish.
Before $\tau$, $E_0\le K\delta^2$,
$\nu\int_0^{t\wedge\tau}G_0\dd r\le K\delta^2$, and $E\le K\alpha^2$.
The thresholds involve continuous adapted processes. On $[0,T]$,
the joint energy controls the physical $H^s$ norm with factor
$(1+T+\rho\sup_{t\le T}|W_t|)^s$, so the continuation criterion excludes
blow-up before $\tau$. After $\tau$, we switch off the nonlinear fluxes
and continue the solution by the homogeneous linear propagators from its
value at $\tau$. These propagators are constructed pathwise from the
continuous frequency characteristics, so the continuation is adapted.
Its zero-mode reference energy is nonincreasing after exit, and its fast
energy satisfies \eqref{eq:critical-linear-drift}. In particular, the
compensated inequality below continues to hold with $a_t=0$ after exit.
Every nonlinear flux is multiplied by $\mathbf1_{t<\tau}$; it is adapted,
and the bounds used below remain valid for all times because the flux
vanishes after exit. The resulting global process agrees with the local
solution before exit, and we retain the same energy notation.

We first combine \eqref{eq:critical-linear-drift} with the refined fast self-interaction bound.
Before exit, $C\nu^{-3}E^2Z\le C\nu^{-3}K^2\alpha^4Z\le Z/4$
by \eqref{eq:split-stopping-smallness}. Except for the term containing
the physical half derivative of the zero field, the remaining coefficients of $E$ are
bounded by $C_{s,\beta}\left[
 \frac{K(\alpha^2+\delta^2)}{\nu}
       +\frac{K^2\delta^4}{\nu^3}\right]$.
By \eqref{eq:split-stopping-smallness}, this coefficient bound is a small multiple
of $\lambda$, since $\nu^{1/2}\le\nu^{1/3}$.
We apply \eqref{eq:separated-mixed-young} to the mixed term. After
absorbing the linear noise drift and dropping the remaining $3Z/4$, we obtain
\eqref{eq:compensated-fast-drift} with
$a_t=C_{s,\beta}\nu^{-1}\mathbf1_{t<\tau}H_0(t)^2$.
The zero-mode energy bounds and \eqref{eq:separated-zero-half} imply
\begin{align*}
 C_t=\int_0^t a_r\dd r
 &\le C_{s,\beta}K\delta^2\nu^{-3/2}\sqrt{t\wedge\tau}\\
 &\le\frac\lambda8t+\Gamma, \quad
 \Gamma=\frac{C_{s,\beta}K^2\delta^4}{\nu^3\lambda}\le1.
\end{align*}
Here $a_t=0$ after exit. The accumulated bound follows from
Cauchy--Schwarz and completion of the square in $\sqrt{t\wedge\tau}$.
Lemma~\ref{lem:compensated-fast} yields
\begin{equation}\label{eq:large-scale-fast-moments}
 (\E E(t)^q)^{1/q}\le C\alpha^2e^{-\lambda t/8},
                          \quad q=1,2,4,
\end{equation}
and
\begin{equation}\label{eq:large-scale-fast-expectation}
 \E\sup_{t\ge0}e^{\lambda t/8}E(t)
       +\E\int_0^\infty e^{\lambda t/8}(\nu D+L)\dd t
                                  \le C_{s,\beta}\alpha^2.
\end{equation}
Lemma~\ref{lem:compensated-fast} also gives
$\E\left[\left(\sup_{t\ge0}e^{\lambda t/8}E(t)\right)^2\right]
                              \le C_{s,\beta}\alpha^4$.
Thus an exit through the fast threshold has probability at
most $C_{s,\beta}/K^2$.

We next use these fast moments to control the simple-zero forcing.
After removing the mean--mean flux, we have
$\norm{\mathcal H_{\rm s}}_{H^s}
               \le C_{s,\beta}(\sqrt K\delta S+E)$.
Mixed zero--fast products have no zero output, and the stopped flux at
$r$ is $\cF_r$-measurable. We can therefore apply conditional divergence
smoothing to this flux. Together with the tower property and Minkowski's
inequality, this gives, for $p=2,4$ and
$x_p(t)=(\E S(t)^p)^{1/p}$,
\begin{align*}
 x_p(t)\le{}&C\eps_0e^{-c_p\mu t}
   +C\sqrt K\delta\int_0^t k_p(t-r)x_p(r)\dd r\\
 &+C\alpha^2\int_0^t k_p(t-r)e^{-\lambda r/8}\dd r,
  \quad k_p(v)=\frac{C_p}{\sqrt{\nu v}}e^{-c_p\mu v}.
\end{align*}
For $p=4$, we use the fourth fast-energy moment.
We choose \[
 0<b_0<\min\{c_2/2,c_4/2,1/16\}.
\]
The feedback kernel then has weighted $L^1$ norm at most $C/\sqrt{\nu\mu}$.
Lemma~\ref{lem:fast-source-convolution} bounds the weighted fast source
by $C\alpha^2/\sqrt{\nu\lambda}$. In view of
\eqref{eq:split-stopping-smallness}, we can absorb the feedback in the
weighted supremum estimate on each finite time interval. Combining this
with the fast-source bound, we obtain
\begin{equation}\label{eq:large-scale-simple-moments}
 x_p(t)\le C_{s,\beta}\delta e^{-b_0\mu t}, \quad p=2,4.
\end{equation}
We now return to the full zero-field energy identity, which has no
martingale in the reference coordinates. Pairing its stopped flux with
$\nabla_tY_0^\sh$ and applying Young's inequality, we obtain
$\dd E_0+\nu G_0\dd t
 \le C_{s,\beta}\nu^{-1}[K\delta^2 S^2+E^2]\dd t$.
Integrating gives, pathwise,
\[
 B_0:=\sup_{t\ge0}E_0(t)+\nu\int_0^\infty G_0\dd t
 \le2\eps_0^2+C_{s,\beta}\nu^{-1}
                     \int_0^\infty(K\delta^2S^2+E^2)\dd t.
\]
To control the probability of exit through the zero threshold, we
estimate the second moment of this pathwise bound. Minkowski's
inequality and the fourth moments in
\eqref{eq:large-scale-fast-moments} and
\eqref{eq:large-scale-simple-moments} give
\begin{align*}
 \left(\E\left[\left(\int_0^\infty S^2\dd t\right)^2\right]\right)^{1/2}
       &\le\int_0^\infty(\E S^4)^{1/2}\dd t\le C\delta^2/\mu,\\
 \left(\E\left[\left(\int_0^\infty E^2\dd t\right)^2\right]\right)^{1/2}
       &\le\int_0^\infty(\E E^4)^{1/2}\dd t\le C\alpha^4/\lambda.
\end{align*}
Combining these moment bounds with the pathwise energy estimate gives
\begin{equation*}
 (\E B_0^2)^{1/2}
 \le C_{s,\beta}\left[
 \eps_0^2+\frac{K\delta^4}{\nu\mu}
                  +\frac{\alpha^4}{\nu\lambda}\right]
 \le C_{s,\beta}\delta^2.
\end{equation*}
Here \eqref{eq:split-stopping-smallness} makes $K\delta^2/(\nu\mu)$
small, while $\alpha^4/(\nu\lambda)\le\delta^2$. This second-moment bound is uniform in
$K$, so the probability of exit through the zero threshold is at most $C_{s,\beta}/K^2$.
The constants in these moment and exit bounds are independent of $K$
under \eqref{eq:split-stopping-smallness}. Choose a fixed
$K_0=K_0(s,\beta)$ large enough and set $K=K_0\zeta^{-1/2}$; the two
exit probabilities then sum to at most $\zeta/4$. If $c$ denotes the
initial-data constant in \eqref{eq:full3d-small}, its assumptions give
\[
 \sqrt K\alpha\le\sqrt{K_0}c\nu^{3/4}, \quad
 \sqrt K\delta\le\sqrt{K_0}(c+c^2)\sigma_0.
\]
After $K_0$ is fixed, choose $c$ so that
$\sqrt{K_0}(c+c^2)\le c_*$. This verifies
\eqref{eq:split-stopping-smallness} without using the desired stability
bounds. On $\{\tau=\infty\}$, the original solution is global. 

With the exit probabilities controlled, we next place the decay and
mean bounds on a common event. Let $0<b<\min\{2b_0,1/8\}$.
The simple-zero energy has no martingale in the reference coordinates.
Using the preceding flux bound, we obtain
\[
 \dd(e^{b\mu t}S^2)+\nu e^{b\mu t}G_{\rm s}\dd t
 \le e^{b\mu t}\left[b\mu S^2
           +C_{s,\beta}\nu^{-1}(K\delta^2S^2+E^2)\right]\dd t.
\]
By \eqref{eq:large-scale-simple-moments} and the $q=2$ case of
\eqref{eq:large-scale-fast-moments},
$\E S(t)^2\le C\delta^2e^{-2b_0\mu t}$ and
$\E E(t)^2\le C\alpha^4e^{-\lambda t/4}$.
Since $b<2b_0$, $b<1/8$ and $\mu\le\lambda$, the weighted source is
integrable. Integrating and taking the time supremum give
\begin{align*}
 \E\left[\sup_{t\ge0}e^{b\mu t}S(t)^2
      +\nu\int_0^\infty e^{b\mu t}G_{\rm s}(t)\dd t\right]
 &\le C_{s,\beta}\left[
 \eps_0^2+\frac{K\delta^4}{\nu\mu}
       +\frac{\alpha^4}{\nu\lambda}+\delta^2\right]
 \le C_{s,\beta}\delta^2.
\end{align*}
The nonlinear mean response $m_{\rm nl}$ solves the heat equation with zero
initial value and the stopped divergence flux. Since $u^2_{00}=0$,
mean--mean products vanish, and mean--oscillatory products average to zero.
The remaining flux has $H^s$ norm at most $C_{s,\beta}(S^2+E)$.
Applying the heat-energy estimate to this flux gives
\begin{equation*}
 \E\left[\sup_{t\ge0}\norm{m_{\rm nl}(t)}_{H^s}^2
   +\nu\int_0^\infty\norm{\partial_y m_{\rm nl}}_{H^s}^2\dd t\right]
 \le C_{s,\beta}\left[
          \frac{\delta^4}{\nu\mu}+\frac{\alpha^4}{\nu\lambda}\right].
\end{equation*}
We apply Markov's inequality to these bounds and \eqref{eq:large-scale-fast-expectation},
allocate $\zeta/12$ to each bound, and intersect the resulting events with
$\{\tau=\infty\}$. On this intersection, the zero threshold gives
\eqref{eq:full3d-zero-budget}, while $\delta\le C\eps$ gives the remaining
energy bounds on an event of probability at least $1-\zeta/2$.
The free heat term contributes at most
$\frac32\|m_{\rm in}\|_{H^s}^2$ to the supremum plus viscous dissipation,
which gives the total mean bound.

Since $\rho/\sqrt\lambda\le C_{s,\beta}$,
Lemma~\ref{lem:physical-fast-damping} supplies a Brownian event with
failure probability at most $\zeta/2$. Intersecting the two events
gives all the stated bounds with probability at least $1-\zeta$. On this
intersection, $B_t\le B_t^3\le C_{s,\beta,\zeta}e^{\lambda t/16}$.
Combining \eqref{eq:physical-fast-frequency-bound} with
\eqref{eq:full3d-fast} and taking square roots gives
\eqref{eq:full3d-physical-damping}.

It remains to prove convergence to zero. Let $\mathcal H_{00}$ denote
the stopped double-mean flux in the heat equation for $m_{\rm nl}$.
The bound $\|\mathcal H_{00}\|_{H^s}\le C_{s,\beta}(S^2+E)$ and the
fourth moments established above give
\[
 \E\int_0^\infty\|\mathcal H_{00}(r)\|_{H^s}^2\dd r
 \le C_{s,\beta}\left(\frac{\delta^4}{\mu}
                         +\frac{\alpha^4}{\lambda}\right)<\infty.
\]
Thus $\mathcal H_{00}\in L^2_tH^s$ almost surely. For $t>T$, the
heat-energy estimate bounds its Duhamel contribution from $[T,t]$ by
\[
 C\nu^{-1/2}
       \|\mathcal H_{00}\|_{L^2([T,\infty);H^s)},
\]
uniformly in $t$. For fixed $T$ and $t\ge2T$, the contribution from $[0,T]$
is bounded by
\[
 \frac{C}{\sqrt{\nu(t-T)}}
       \|\mathcal H_{00}\|_{L^1([0,T];H^s)}\longrightarrow0.
\]
The initial heat term tends strongly to zero by Fourier dominated
convergence on $\R_y$. Taking $t\to\infty$ first and then $T\to\infty$
proves convergence of the double mean. The remaining modes decay
exponentially by \eqref{eq:full3d-fast} and \eqref{eq:full3d-simple}, so the
full reference field converges in $H^s$.

The double mean is unchanged by the affine pullback and therefore also
converges in physical $H^s$. For the nonmean part, the comparison
\[
 \|f_{\rm osc}(t)\|_{H^s}
 \le C_s(1+t+\rho|W_t|)^s\|f_{\rm osc}^\sh(t)\|_{H^s},
  \quad f_{\rm osc}=f-f_{00},
\]
multiplies only an exponentially decaying norm. Since Brownian paths
grow subexponentially almost surely, this physical norm tends to zero.
The separate scales follow from \eqref{eq:full3d-small}.
\end{proof}

\Needspace{6\baselineskip}
\subsection{The temperature mean in three dimensions}\label{sec:temperature-mean}

We finally apply the thermal normal form to the full system. Let a bar denote
normalized $x,z$ averaging. Every local strong solution satisfies
\begin{equation}\label{eq:full-temperature-normal-form}
 \begin{split}
 (\partial_t-\nu\partial_{yy})
   \left(\overline\theta+
          \frac{1}{2\beta}\partial_y\overline{\theta^2}\right)
 ={}&-\frac{1}{2\beta}\partial_{yy}\overline{u^2\theta^2} -\frac{\nu}{\beta}\partial_y\overline{|\nabla\theta|^2}.
 \end{split}
\end{equation}
Indeed, the product calculation in
Lemma~\ref{lem-scalar-normal-form} applies, since averaging removes both affine
scalar transports and $\overline{u^2}=0$. The streamwise pressure
derivative prevents the same extension of the $q$ identity.

\Needspace{6\baselineskip}

\begin{proof}[Proof of Corollary~\ref{cor:full-temperature-mean}]
We use the main theorem at index $s+2$ to obtain
\[
 \begin{gathered}
 \sup_t\|U^\sh\|_{H^{s+2}}\le C\eps, \quad
 \nu\int_0^\infty\|\nabla_tU^\sh\|_{H^{s+1}}^2\dd t\le C\eps^2,\\
 \int_0^\infty\|U^\sh-\overline U\|_{H^{s+2}}^2\dd t
       \le C(\delta^2/\mu+\alpha^2/\lambda).
 \end{gathered}
\]
Averaging commutes with the affine pullback. Since $\overline{u^2}=0$,
expanding $\theta$ into mean and oscillatory parts shows that each term
in $\overline{u^2\theta^2}$ has at least two oscillatory factors.
Its $L^1_tH^{s+2}$ norm is therefore at most
$C\eps(\delta^2/\mu+\alpha^2/\lambda)\le C\eps^2$:
the data conditions imply $\delta\le C\eps$,
$\delta\lesssim\sqrt{\nu\mu}\le\mu$, and
$\alpha\lesssim\nu^{3/4}\le\lambda$.
The boundary corrections in \eqref{eq:full-temperature-normal-form}
have $H^s$ norm at most $C\eps^2$, and the viscosity-weighted gradient
term has $L^1_tH^s$ norm at most $C\eps^2$.
Subtracting $e^{\nu t\partial_{yy}}\overline{\theta_{\rm in}}$ from the
heat Duhamel formula, we use contractivity of the heat semigroup to bound
the remaining initial quadratic correction. Combining these bounds gives
\eqref{eq:full-temperature-mean}.
\end{proof}

\section*{Acknowledgements}
This work was partially supported by the National Key Research and
Development Program of China (Grant No.~2023YFC2206100), the National
Natural Science Foundation of China (Grant No.~12231008), and the
Fundamental Research Funds for the Central Universities
(Grant No.~2026BRSXB003).

\section*{Conflict of interest}
The author declares no conflict of interest.

\section*{Data availability}
No data were generated or used in this study.


\end{document}